\documentclass[12pt]{amsart}%
\usepackage{amsmath}
\usepackage{graphicx}
\usepackage[small,nohug,heads=littlevee]{diagrams}
\usepackage{amscd}
\usepackage{geometry}
\usepackage{amsfonts}
\usepackage{amssymb}%
\diagramstyle[labelstyle=\scriptstyle]
\newtheorem{theorem}{Theorem}[section]
\theoremstyle{plain}

\newtheorem{corollary}[theorem]{Corollary}

\newtheorem{example}{Example}

\newtheorem{lemma}[theorem]{Lemma}

\newtheorem{proposition}[theorem]{Proposition}
\newtheorem{remark}{Remark}

\numberwithin{equation}{section}

\begin{document}
\title[Topological properties of spaces admitting free group actions]{Topological properties of spaces admitting free group actions}
\author{Ross Geoghegan}
\address{Department of Mathematical Sciences, Binghamton University, Binghamton, New
York 13902-6000}
\email{ross@math.binghamton.edu}
\author{Craig R. Guilbault}
\address{Department of Mathematical Sciences, University of Wisconsin-Milwaukee,
Milwaukee, Wisconsin 53201}
\email{craigg@uwm.edu}
\thanks{This project was aided by a Simons Foundation Collaboration Grant awarded to
the second author. }
\date{October 15, 2011}

\begin{abstract}
In 1992, David Wright proved a remarkable theorem about which contractible
open manifolds are covering spaces. He showed that if a one-ended open
manifold $M^{n}$ has pro-mono\-mor\-phic fundamental group at infinity which
is not pro-trivial and is not stably ${\mathbb{Z}}$, then $M$ does not cover
any manifold (except itself). In the non-manifold case, Wright's method showed
that when a one-ended, simply connected, locally compact ANR $X$ with
pro-mono\-mor\-phic fundamental group at infinity admits an action of
$\mathbb{Z}$ by covering transformations then the fundamental group at
infinity of $X$ is (up to pro-isomorphism) an inverse sequence of finitely
generated free groups. We improve upon this latter result, by showing that $X$
must have a \emph{stable} finitely generated free fundamental group at
infinity. Simple examples show that a free group of any finite rank is possible.

We also prove that if $\ X$ (as above), admits a non-cocompact action of
$\mathbb{Z\times Z}$ by covering transformations, then $X$ is simply connected
at infinity.

We deduce the following corollary in group theory: Every finitely presented
one-ended group $G$ which contains an element of infinite order satisfies
exactly one of the following:

\begin{itemize}
\item $G$ is simply connected at infinity;

\item $G$ is virtually a surface group;

\item The fundamental group at infinity of $G$ is not pro-mono\-mor\-phic.

\end{itemize}

Our methods also provide a quick new proof of Wright's open manifold theorem.

\end{abstract}
\maketitle

\section{Introduction}

In this paper we address a series of questions which are of interest in both
topology and geometric group theory.

\medskip

Let $X$ be a locally finite simply connected one-ended CW complex (or, more
generally, a locally compact simply connected one-ended metric ANR). We are
interested in the fundamental group at infinity\footnote{This and several
other concepts mentioned in the introduction will be discussed in detail later
in the paper.} of $X$, and what special properties it must have in order that
$X$ can be a non-trivial covering space. In particular, we consider:

\medskip

\noindent\textbf{Question 1:} \emph{Does }$\mathbb{Z}$ \emph{act as a group of
covering transformations on }$X$?

\medskip

\noindent\textbf{Question 2:} \emph{Does some group }$G$ \emph{act cocompactly
as a group of covering transformations on }$X$?

\medskip

Turning the latter question around, let $K$ be a compact ANR (for example, a
finite complex) with one-ended universal cover $X$. Can the fundamental group
at infinity of $X$ be arbitrary? Or how restricted is it? \medskip

Here are closely related manifold versions of our questions: \medskip

\noindent\textbf{Question 3:} \emph{If} $M^{n}$ \emph{is a contractible open
manifold which is not homeomorphic to} ${\mathbb{R}}^{n}$ \emph{does it cover
any manifold non-trivially? In other words, does} $M^{n}$ \emph{support a
properly discontinuous free action of a non-trivial (necessarily torsion free)
group?}

\medskip\noindent\textbf{Question 4:} \emph{Must the universal cover of a
closed aspherical} $n$\emph{-manifold be homeomorphic to} $\mathbb{R}^{n}$?

\medskip

We are certainly not the first authors to address these questions. We will
begin by reviewing some of what is known, starting with the manifold case.
\medskip

\noindent\textbf{Concerning Question 4}

\medskip

In dimensions $\leq3$ the universal cover of a closed aspherical manifold is
indeed homeomorphic to a Euclidean space. This is classical in dimensions
$\leq2$ and follows from Perelman's solution to the Poincar\'{e} Conjecture in
dimension 3. A negative answer was obtained by Davis \cite{Da} in all
dimensions $\geq4$. The invariant which detects Davis's examples is the
fundamental group at infinity; specifically, Davis gave examples of a
one-ended Coxeter group $\Gamma$ having a subgroup $G$ of finite index such
that the universal cover of a closed manifold $K(G,1)$ is not simply connected
at infinity.

\medskip

\noindent\textbf{Concerning Question 3}

\medskip

Recall that when $n\geq3$ a contractible open $n$-manifold is homeomorphic to
$\mathbb{R}^{n}$ if and only if it is simply connected at infinity. (This will
be discussed further in \S \ref{Sec: Manifold results}.) So one looks for an
answer among contractible manifolds which are not simply connected at infinity
(i.e. with non-trivial fundamental group at infinity). Whitehead \cite{Wh}
gave the first example of a contractible open $3$-manifold which is not simply
connected at infinity. Later, it was shown in \cite{Mc} that uncountably many
pairwise non-homeomorphic such $3$-manifolds exist. But until the 1980's it
was not known whether any of these could cover a manifold non-trivially.
\medskip

Myers \cite{My1} obtained the first notable result along these lines, proving
that no member of a certain class of Whitehead-type contractible open
$3$-manifolds admits an action of $\mathbb{Z}$ by covering transformations.
Later, Wright \cite{Wr} gave a significant generalization by showing that a
contractible open $n$-manifold ($n\geq3$) with pro-\allowbreak
mono\-mor\-phic\footnote{This is defined in \S \ref{Sec: background}; roughly,
it means that the inverse sequence of fundamental groups of complements of
larger and larger compact submanifolds looks like a sequence of
monomorphisms.} fundamental group at infinity, which admits a nontrivial
action by covering transformations, is necessarily simply connected at
infinity. All the aforementioned Whitehead-type $3$-manifolds have this
pro-\allowbreak mono\-mor\-phic property, hence they do not cover manifolds
non-trivially. Wright's theorem also implies that the interiors of compact
contractible $n$-manifolds with non-simply connected boundaries do not cover
non-trivially. \medskip

\noindent\textbf{Concerning Question 1} \medskip

Wright's method extends beyond manifold topology\footnote{Background and more
detail concerning Wright's Theorem can be found in Section 16.3 of \cite{Ge}%
.}. In particular, he proved:

\begin{theorem}
[{\cite[Th.9.1]{Wr}}]\label{Th: Wright's pro-free theorem} Let $X$ be a simply
connected one-ended locally compact ANR, and assume $X$ has pro-\allowbreak
mono\-mor\-phic fundamental group at infinity. If $\mathbb{Z}$ acts as
covering transformations on $X$ then the fundamental group at infinity of $X$
is pro-free and pro-finitely generated.
\end{theorem}

\medskip

For this version of Wright's Theorem, see {\cite[Th.16.3.4]{Ge}}.

\medskip

The conclusion of Theorem \ref{Th: Wright's pro-free theorem} implies that the
fundamental group at infinity may be represented by
\[
F_{1}\leftarrow F_{2}\longleftarrow F_{3}\leftarrow\cdots
\]
where each $F_{i}$ is a finitely generated free group. Given the
pro-\allowbreak mono\-mor\-phic hypothesis, one may also assume that the
bonding homomorphisms in this sequence are injective. The main theorem of this
paper is an improvement of Wright's Theorem
\ref{Th: Wright's pro-free theorem}, namely:

\begin{theorem}
\label{Th: Stably free theorem}Let $X$ be a simply connected one-ended locally
compact ANR with pro-\allowbreak mono\-mor\-phic fundamental group at
infinity. If $\mathbb{Z}$ acts as covering transformations on $X$, then the
fundamental group at infinity of $X$ is stably a finitely generated free group.
\end{theorem}

The conclusion of Theorem \ref{Th: Stably free theorem} means that the
finitely generated free groups $F_{i}$ all can be taken to have the same rank,
and the bonding morphisms can be taken to be isomorphisms. \medskip

\begin{remark}
\emph{\textquotedblleft Stable\textquotedblright\ means, roughly, that the
inverse sequence of fundamental groups of complements of larger and larger
compacta looks like a sequence of isomorphisms. Thus \textquotedblleft
stable\textquotedblright\ is equivalent to \textquotedblleft pro-monomorphic
and pro-epimorphic\textquotedblright. \textquotedblleft
Pro-epimorphic\textquotedblright\ is also known as \textquotedblleft
semistable\textquotedblright\ or \textquotedblleft
Mittag-Leffler\textquotedblright. Thus Theorem \ref{Th: Stably free theorem}
improves on Wright's Theorem by establishing semistability.}
\end{remark}

\medskip

\begin{example}
As a simple illustration, let $X$ be a simply connected polyhedron with
``solenoidal" fundamental group at infinity; i.e. having the form
\[
\mathbb{Z}\overset{\times2}{\longleftarrow}\mathbb{Z}\overset{\times
2}{\longleftarrow}\mathbb{Z}\overset{\times2}{\longleftarrow}\cdots.
\]
(Such spaces are easy to construct.) Our theorem prohibits $X$ from admitting
an action of $\mathbb{Z}$ by covering transformations.
\end{example}

\medskip

On the other hand, there exist contractible spaces which have stable free
fundamental groups at infinity of any finite rank, and which admit actions of
$\mathbb{Z}$ by covering transformations: \medskip

\begin{example}
Let $K_{n}$ be a wedge of $n$ rays with a common vertex and let $X=K_{n}%
\times\mathbb{R}$. Then $X$ has a stable fundamental group at infinity that is
free of rank $n-1$. Moreover, the homeomorphism $j:X\rightarrow X$ induced by
a translation in the $\mathbb{R}$-factor generates an action of $\mathbb{Z}$
by covering transformations. As a variation, one may combine translation in
the $\mathbb{R}$-coordinate with a permutation of the rays of $K_{n}$. By
doing so, it is possible to obtain quotient spaces with different numbers of
ends and various fundamental group behavior at those ends. These examples will
be useful to keep in mind when reading
\S \ref{Section: Rectangular neighborhoods}.%
%TCIMACRO{\FRAME{ftbpFU}{4.7721in}{0.8997in}{0pt}{\Qcb{$K_{3}\times\mathbb{R}%
%.$}}{\Qlb{fig1}}{actions-fig1.eps}{\special{ language "Scientific Word";
%type "GRAPHIC";  maintain-aspect-ratio TRUE;  display "USEDEF";
%valid_file "F";  width 4.7721in;  height 0.8997in;  depth 0pt;
%original-width 7.9083in;  original-height 1.4529in;  cropleft "0";
%filename 'actions-fig1.eps';file-properties "XNPEU";}}}%
%BeginExpansion
\begin{figure}
[ptb]
\begin{center}
\includegraphics[
height=0.8997in,
width=4.7721in
]%
{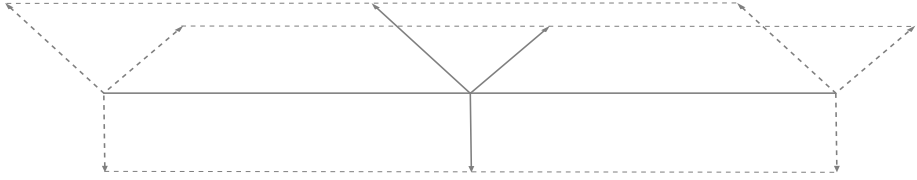}%
\caption{$K_{3}\times\mathbb{R}.$}%
\label{fig1}%
\end{center}
\end{figure}
%EndExpansion

\end{example}

\medskip

When $X$ admits a free $(\mathbb{Z\times Z)}$-action, we are able to prove a
stronger result.

\begin{theorem}
\label{Th: Z+Z theorem}Let $X$ be a simply connected one-ended locally compact
ANR with pro-\allowbreak mono\-mor\-phic fundamental group at infinity. If $X$
admits an action of $\mathbb{Z\times Z}$ by covering transformations, then
either that action is cocompact or $X$ is simply connected at infinity.
\end{theorem}

\medskip

\noindent\textbf{Concerning Question 2}

\medskip

Let $G$ be a finitely presented one-ended group having an element of infinite
order, and let $K$ be a finite complex whose fundamental group is isomorphic
to $G$. The universal cover $X$ of $K$ has fundamental group at infinity
represented by an inverse sequence%

\[
G_{1}\longleftarrow G_{2}\longleftarrow G_{3}\longleftarrow\cdots
\]
where each $G_{i}$ is the fundamental group of the complement of a compact
subcomplex of $X$. No example is known where this sequence fails to be
semistable. Easy examples occur where this sequence (up to pro-isomorphism)
has or does not have pro-\allowbreak mono\-mor\-phic fundamental group at
infinity. Here, we only consider the case where $X$ has pro-\allowbreak
mono\-mor\-phic fundamental group at infinity. Theorem
\ref{Th: Stably free theorem} implies that, letting $j$ be an element of
infinite order, the infinite cyclic group $\left\langle j\right\rangle $ acts
as covering transformations on $X$, hence $X$ has a stable finitely generated
free fundamental group at infinity. In other words, one may assume the groups
$G_{i}$ are all free and finitely generated of the same rank and the bonding
morphisms are isomorphisms. But in this case more is known. As explained in
the proof of {\cite[Th.16.5.6]{Ge}}, for homological reasons a theorem of
Farrell \cite{Fa} implies that this rank must be either 0 or 1. Moreover, in
the rank 1 case, Bowditch \cite{Bo} has shown that $G$ is virtually a surface
group---meaning that $G$ contains a finite index subgroup which is the
fundamental group of a closed surface. In other words, a consequence of our
main theorem is:

\begin{theorem}
\label{cocompact} Let the one-ended finitely presented group $G$ have
pro-\allowbreak mono\-mor\-phic fundamental group at infinity and assume $G$
contains an element of infinite order. Then $G$ is either simply connected at
infinity or $G$ is virtually a surface group.
\end{theorem}

\begin{remark}
\emph{The advance, here, over \cite[Th.16.5.6]{Ge} is that we do not have to
make semistability of }$G$\emph{ part of the hypothesis.}
\end{remark}

\begin{remark}
\emph{The one-ended examples of Davis, mentioned above, where }$G$\emph{ is a
torsion free subgroup of finite index in a suitable Coxeter group, show that
the universal cover }$X$\emph{ in Theorem \ref{cocompact} does not always have
pro-\allowbreak mono\-mor\-phic fundamental group at infinity when it is the
universal cover of a finite complex.}
\end{remark}

\medskip

A similar application of Theorem \ref{Th: Stably free theorem} applies to
CAT(0) groups.

\begin{theorem}
If the CAT(0) group $G$ acts geometrically on a one-ended proper CAT(0) space
$X$ such that $\partial_{\infty}X$ is 1-dimensional and $X$ has
pro-\allowbreak mono\-mor\-phic fundamental group at infinity, then $G$ is
virtually a surface group.
\end{theorem}

\begin{proof}
\cite[Th. 11]{Sw} ensures the existence of an infinite order element in $G$,
thus permitting an application of Theorem \ref{Th: Stably free theorem}. Then
\cite[Main Theorem]{GO} rules out the possibility of $X$ being simply
connected at infinity when $\partial_{\infty}X$ is 1-dimensional.
\end{proof}

\medskip

This paper also includes some manifold theoretic results, most notably, a
brief new proof of \cite[Main Theorem]{Wr} and a new proof of a theorem from
\cite{GM}.

\medskip

\noindent\textbf{A special case of the Borel Construction}

\medskip

The starting point of our argument is the following useful observation:

\begin{proposition}
\label{Prop: mapping torus observation}Let $Y$ be a connected, locally path
connected space on which an infinite cyclic group $J=\left\langle
j\right\rangle $ acts as covering transformations. Then $(J\backslash
X)\times\mathbb{R}$ is homeomorphic to the mapping torus $T_{j}\left(
X\right)  $.
\end{proposition}

\begin{remark}
\emph{In \cite{GM}, Proposition \ref{Prop: mapping torus observation} is
proved using bundle theory; the authors cite an earlier paper by Farrell
\cite{Fa}. Another proof is given in \cite{Ge} as Proposition 13.7.4. A new
and elementary proof of Proposition \ref{Prop: mapping torus observation} is
contained within our proof of Theorem \ref{Th: Stably free theorem}; see
Remark \ref{Remark: elementary proof} in
\S \ref{Sec: Proof of Theorem-topological part}.}
\end{remark}

Proposition \ref{Prop: mapping torus observation} provides two pictures of the
same space. Those pictures yield canonical, but very different, families of
neighborhoods of infinity---\emph{rectangular} neighborhoods of infinity in
the case of $(J\backslash X)\times\mathbb{R}$ and \emph{mapping torus}
neighborhoods of infinity in the case of $T_{j}\left(  X\right)  $ (see
\S \ref{Section: Rectangular neighborhoods} and
\S \ref{Section: Mapping torus neighborhoods}). By comparing the fundamental
group systems arising from these different pictures, we are sometimes able to
coax out information about the original space $X$, as was done, for example,
in \cite{GM}. The delicate nature of that task accounts for much of the work
found here.

\medskip\noindent\textbf{Layout} \medskip

The layout of this paper is as follows. In \S \ref{Sec: background} we provide
most of the necessary definitions and background. The first portion of that
section is purely algebraic---dealing primarily with inverse sequences of
groups. In addition to essential definitions and notation, some basic results
are proved for later use. The latter portion of \S \ref{Sec: background}
discusses the topology of noncompact spaces; neighborhoods of infinity and
fundamental pro-groups are discussed and some useful equivalences are
reviewed. In \S \ref{Section: Rectangular neighborhoods} and
\S \ref{Section: Mapping torus neighborhoods} we look at canonical
neighborhoods of infinity in the homeomorphic spaces, $(J\backslash
X)\times\mathbb{R}$ and $T_{j}\left(  X\right)  $, of Proposition
\ref{Prop: mapping torus observation}. Of particular importance will be some
clean descriptions of the fundamental groups of those neighborhoods of
infinity. The remainder of the paper contain the proofs of our principal
results. In \S \ref{Sec: Manifold results} we focus on manifolds. Our new
proof of Wright's main theorem may be viewed as a warm-up for the more general
(non-manifold) version that is Theorem \ref{Th: Stably free theorem}.
\S \ref{Sec: Proof of Theorem-algebraic part} and
\S \ref{Sec: Proof of Theorem-topological part} are devoted to the proof of
Theorem \ref{Th: Stably free theorem}. In the final section we push our
techniques one step further to prove Theorem \ref{Th: Z+Z theorem}.

The authors wish to acknowledge the contribution of the referee, whose
comments led to a significant simplification in the proof of our main theorem.

\section{Definitions and Background\label{Sec: background}}

This section contains much of the terminology and notation needed for the
remainder of the paper. In addition, several preliminary results are
presented. The section is divided into two parts; the first is entirely
algebraic---dealing primarily with inverse sequences of groups, while the
second is topological---dealing primarily with the topology at the ends of
noncompact spaces.

\subsection{Algebra of inverse sequences}

Throughout this subsection all arrows denote homomorphisms, while arrows of
the type $\twoheadrightarrow$ or $\twoheadleftarrow$ denote surjections and
arrows of the type $\rightarrowtail$ and $\leftarrowtail$ denote injections.
The symbol $\cong$ indicates an isomorphism.

Let
\[
G_{0}\overset{\lambda_{1}}{\longleftarrow}G_{1}\overset{\lambda_{2}%
}{\longleftarrow}G_{2}\overset{\lambda_{3}}{\longleftarrow}\cdots
\]
be an inverse sequence of groups. A \emph{subsequence} of $\left\{
G_{i},\lambda_{i}\right\}  $ is an inverse sequence of the form
\[
\begin{diagram}
G_{i_{0}} & \lTo^{\lambda_{i_{0}+1}\circ\cdots\circ\lambda_{i_{1}}
} & G_{i_{1}} & \lTo^{\lambda_{i_{1}+1}\circ\cdots\circ
\lambda_{i_{2}}} & G_{i_{2}} & \lTo^{\lambda_{i_{2}+1}\circ
\cdots\circ\lambda_{i_{3}}} & \cdots.
\end{diagram}
\]
In the future we will denote a composition $\lambda_{i}\circ\cdots\circ
\lambda_{j}$ ($i\leq j$) by $\lambda_{i,j}$.

Sequences $\left\{  G_{i},\lambda_{i}\right\}  $ and $\left\{  H_{i},\mu
_{i}\right\}  $ are \emph{pro-isomorphic} if, after passing to subsequences,
there exists a commuting \textquotedblleft ladder diagram\textquotedblright:
\begin{equation}
\begin{diagram} G_{i_{0}} & & \lTo^{\lambda_{i_{0}+1,i_{1}}} & & G_{i_{1}} & & \lTo^{\lambda_{i_{1}+1,i_{2}}} & & G_{i_{2}} & & \lTo^{\lambda_{i_{2}+1,i_{3}}}& & G_{i_{3}}& \cdots\\ & \luTo & & \ldTo & & \luTo & & \ldTo & & \luTo & & \ldTo &\\ & & H_{j_{0}} & & \lTo^{\mu_{j_{0}+1,j_{1}}} & & H_{j_{1}} & & \lTo^{\mu_{j_{1}+1,j_{2}}}& & H_{j_{2}} & & \lTo^{\mu_{j_{2}+1,j_{3}}} & & \cdots \end{diagram} \label{basic ladder diagram}%
\end{equation}
Clearly an inverse sequence is pro-isomorphic to any of its subsequences. To
avoid tedious notation, we sometimes do not distinguish $\left\{
G_{i},\lambda_{i}\right\}  $ from its subsequences. Instead we assume that
$\left\{  G_{i},\lambda_{i}\right\}  $ has the desired properties of a
preferred subsequence---prefaced by the words \textquotedblleft after passing
to a subsequence and relabeling\textquotedblright.

The \emph{inverse limit }of a sequence $\left\{  G_{i},\lambda_{i}\right\}  $
is a subgroup of $\prod G_{i}$ defined by
\[
\underleftarrow{\lim}\left\{  G_{i},\lambda_{i}\right\}  =\left\{  \left.
\left(  g_{0},g_{1},g_{2},\cdots\right)  \in\prod_{i=0}^{\infty}%
G_{i}\right\vert \lambda_{i}\left(  g_{i}\right)  =g_{i-1}\right\}  .
\]
Notice that, for each $i$, there is a \emph{projection homomorphism}
$p_{i}:\underleftarrow{\lim}\left\{  G_{i},\lambda_{i}\right\}  \rightarrow
G_{i}$. It is a standard fact that pro-isomorphic inverse sequences have
isomorphic inverse limits.

An inverse sequence $\left\{  G_{i},\lambda_{i}\right\}  $ is \emph{stable} if
it is pro-isomorphic to a constant inverse sequence $\left\{
H,\operatorname{id}_{H}\right\}  $, or equivalently, to an inverse sequence
$\left\{  H_{i},\mu_{i}\right\}  $ where each $\mu_{i}$ is an isomorphism. In
those cases, the projection homomorphisms take $\underleftarrow{\lim}\left\{
G_{i},\lambda_{i}\right\}  $ isomorphically onto $H$ and each of the$\ H_{i}$.

Another condition equivalent to the stability of $\left\{  G_{i},\lambda
_{i}\right\}  $ is that, after passing to an appropriate subsequence and
relabeling, there exists a commutative diagram of the form
\begin{equation}
\begin{diagram} G_{0}& & \lTo^{{\lambda}_{1}} & & G_{1} & & \lTo^{{\lambda}_{2}} & & G_{2} & & \lTo^{{\lambda}_{3}} & & G_{3} &\cdots\\ & \luTo & & \ldTo & & \luTo & & \ldTo & & \luTo & & \ldTo & \\ & & \operatorname{Im}\left( \lambda_{1}\right) & & \lTo^{\cong} & & \operatorname{Im}\left( \lambda _{2}\right) & &\lTo^{\cong} & & \operatorname{Im}\left( \lambda_{3}\right) & & \lTo^{\cong} & &\cdots & \\ \end{diagram} \label{Standard stability ladder}%
\end{equation}
where all unlabeled homomorphisms are obtained by restriction or inclusion.

If $\left\{  G_{i},\lambda_{i}\right\}  $ is pro-isomorphic to some $\left\{
H_{i},\mu_{i}\right\}  $, where each $\mu_{i}$ is an epimorphism, we call
$\left\{  G_{i},\lambda_{i}\right\}  $ \emph{semi\-stable }(or
\emph{Mittag-Leffler, }or\emph{\ pro-epimorphic}). In that case, there exists
a ladder diagram of the type described in ($\ref{Standard stability ladder}$),
but with the isomorphisms replaced by epimorphisms. Similarly, if $\left\{
H_{i},\mu_{i}\right\}  $ can be chosen so that each $\mu_{i}$ is a
monomorphism, we call our inverse sequence \emph{pro-\allowbreak
mono\-mor\-phic}; in that case, there exists a diagram of type
(\ref{Standard stability ladder}) for which maps in the bottom row are
monomorphisms. It is easy to show that an inverse sequence that is both
semi\-stable and pro-\allowbreak mono\-mor\-phic is stable.

Given a group $G$, a pair of isomorphic subgroups $H,H^{\prime}\leq G$, and a
specified isomorphism $\varphi:H\rightarrow H^{\prime}$, the \emph{HNN
extension} of $G$ by $\varphi$ is the group
\[
G\ast_{\varphi}=\left\langle G,t\mid t^{-1}ht=\varphi(h)\text{ }\forall\text{
}h\in H\right\rangle .
\]
In this setup, $G$ is called the \emph{base group}, $H$ and $H^{\prime}$ the
\emph{associated subgroups}, and $t$ the \emph{stable letter}. It is a
standard fact that $G$ injects into $G\ast_{\varphi}$ (so we view $G$ as a
subgroup), that $t$ generates an infinite cyclic subgroup $\left\langle
t\right\rangle \leq G\ast_{\varphi}$, and that $G\cap\left\langle
t\right\rangle =\left\{  1\right\}  $. We refer to \cite{Co} or \cite{LS} for
these and additional properties of HNN extensions.

A special case of HNN extension occurs when $H=H^{\prime}=G$. In that
situation, the HNN extension is a \emph{semidirect product} of $G$ with
$\left\langle t\right\rangle $ with respect to the automorphism $\varphi
:G\rightarrow G$. That group will be denoted $G\rtimes_{\varphi}\left\langle
t\right\rangle $ or---when the specific isomorphism is not important---just
$G\rtimes\left\langle t\right\rangle $. A nice exposition of semidirect
products is in \cite[Ch.8]{Me}.

Some of the topological constructions used in this paper produce entire
inverse sequences of HNN extensions and/or semidirect products from an initial
\textquotedblleft base\textquotedblright\ sequence. Of special interest here
is the extent to which a given property of one of those sequences implies the
same property for the other. Before formulating some propositions of that
sort, we state a pair of elementary observations and a pair of corollaries
that will be used throughout; proofs are left as exercises.

\begin{lemma}
\label{induced HNN homomorphism}Let $G_{0}$ and $G_{1}$ be groups and
$\varphi_{i}:H_{i}\rightarrow H_{i}^{\prime}$ an isomorphism between subgroups
of $G_{i}$ for $i=0,1$. Suppose $\lambda:G_{1}\rightarrow G_{0}$ is a
homomorphism taking $H_{1}$ into $H_{0}$ and $H_{1}^{\prime}$ into
$H_{0}^{\prime}$ such that $\varphi_{0}\circ\lambda\left(  h\right)
=\lambda\circ\varphi_{1}\left(  h\right)  $ for all $h\in H_{1}$. Then
$\lambda$ induces a unique homomorphism $\bar{\lambda}:G_{1}\ast_{\varphi_{1}%
}\rightarrow G_{0}\ast_{\varphi_{0}}$ that restricts to $\lambda$ on $G_{1}$
and sends the stable letter $t_{1}$ to the stable letter $t_{0}$.
\end{lemma}

\begin{lemma}
\label{Lem: HNN induced triangle}Given the above setup, the homomorphism
$\varphi_{0}$ restricts to an isomorphism $\psi_{0}$ of $\lambda\left(
H_{1}\right)  $ onto $\lambda\left(  H_{1}^{\prime}\right)  $. If we let
$\upsilon:\operatorname{Im}\left(  \lambda\right)  \hookrightarrow G_{0}$ be
inclusion and $\delta:G_{1}\rightarrow\operatorname{Im}\left(  \lambda\right)
$ the corestriction of $\lambda$, then the induced homomorphisms provided by
Lemma \ref{induced HNN homomorphism} yield a commutative diagram:%
\[
\begin{diagram}
G_{0}\ast_{\varphi_{0}} & & \lTo^{\overline{\lambda}}  & &
G_{1}\ast_{\varphi_{1}}\\
& \luTo^{\overline{\upsilon}} & & \ldTo^{\overline{\delta}} & \\
& & \operatorname{Im}\left(  \lambda\right)  \ast_{\psi_{0}} & &\\
\end{diagram}
\]

\end{lemma}

\begin{corollary}
\label{Cor: induced semidirect product homomorphism}If, under the setup of
Lemma \ref{induced HNN homomorphism}, $H_{i}=H_{i}^{\prime}=G_{i}$ for each
$i$, then $\lambda$ induces a unique homomorphism $\bar{\lambda}:G_{1}%
\rtimes_{\varphi_{1}}\left\langle t_{1}\right\rangle \rightarrow G_{0}%
\rtimes_{\varphi_{0}}\left\langle t_{0}\right\rangle $ that restricts to
$\lambda$ on $G_{1}$ and takes $t_{1}$ to $t_{0}$.
\end{corollary}

\begin{corollary}
\label{Cor: Semidirect product-induced triangle}Given the setup of Corollary
\ref{Cor: induced semidirect product homomorphism}, the homomorphism
$\varphi_{0}$ restricts to an automorphism $\psi_{0}$ of $\operatorname{Im}%
\left(  \lambda\right)  $. If we let $\upsilon:\operatorname{Im}\left(
\lambda\right)  \hookrightarrow G_{0}$ be inclusion and $\delta:G_{1}%
\rightarrow\operatorname{Im}\left(  \lambda\right)  $ the corestriction of
$\lambda$, then the induced homomorphisms provided by Corollary
\ref{Cor: induced semidirect product homomorphism} yield a commutative
diagram:%
\[
\begin{diagram}
G_{0}\rtimes_{\varphi_{0}}\left\langle t_{0}\right\rangle & & \lTo^{\overline{\lambda}}  & &
G_{1}\rtimes_{\varphi_{1}}\left\langle t_{1}\right\rangle\\
& \luTo^{\overline{\upsilon}} & & \ldTo^{\overline{\delta}} & \\
& & \operatorname{Im}\left(  \lambda\right)  \rtimes_{\psi_{0}}\left\langle
t_{0}\right\rangle & &\\
\end{diagram}
\]
Since $\overline{\upsilon}$ is injective, we view it as an inclusion map.
\end{corollary}

For the remainder of this section, we work with the following setup:

\begin{enumerate}
\item[(a)] $G_{0}\overset{\lambda_{1}}{\longleftarrow}G_{1}\overset
{\lambda_{2}}{\longleftarrow}G_{2}\overset{\lambda_{3}}{\longleftarrow}\cdots$
is an inverse sequence of groups (the \emph{base} sequence),

\item[(b)] for each $i$, $H_{i},H_{i}^{\prime}\leq G_{i}$ and $\varphi
_{i}:H_{i}\rightarrow H_{i}^{\prime}$ is an isomorphism, and

\item[(c)] for each $i$, $\varphi_{i-1}\circ\lambda_{i}\left(  h\right)
=\lambda_{i}\circ\varphi_{i}\left(  h\right)  $ for all $h\in H_{i}$.
\end{enumerate}

\noindent Under these assumptions, repeated application of Lemma
\ref{induced HNN homomorphism} gives the following \emph{induced HNN
sequence}:%
\[
G_{0}\ast_{\varphi_{0}}\overset{\bar{\lambda}_{1}}{\longleftarrow}G_{1}%
\ast_{\varphi_{1}}\overset{\bar{\lambda}_{2}}{\longleftarrow}G_{2}%
\ast_{\varphi_{2}}\overset{\bar{\lambda}_{3}}{\longleftarrow}\cdots\text{.}%
\]
If we strengthen condition (b) to:

\begin{enumerate}
\item[(b$^{\prime}$)] for each $i$, $\varphi_{i}\in\operatorname{Aut}\left(
G_{i}\right)  $,
\end{enumerate}

\noindent then repeated application of Corollary
\ref{Cor: induced semidirect product homomorphism} gives the following
\emph{induced semidirect product sequence}:%
\begin{equation}
G_{0}\rtimes_{\varphi_{0}}\left\langle t_{0}\right\rangle \overset
{\bar{\lambda}_{1}}{\longleftarrow}G_{1}\rtimes_{\varphi_{1}}\left\langle
t_{1}\right\rangle \overset{\bar{\lambda}_{2}}{\longleftarrow}G_{2}%
\rtimes_{\varphi_{2}}\left\langle t_{2}\right\rangle \overset{\bar{\lambda
}_{3}}{\longleftarrow}\cdots\text{.}%
\end{equation}

The situation for semidirect products is less complicated than for general HNN
extensions, so we turn to that case first.

\begin{proposition}
\label{Prop: semidirect vs base sequences}Assume we are given the above setup,
with Condition (b$^{\prime}$) in place of (b). Then

\begin{enumerate}
\item The induced semidirect product sequence is pro-\allowbreak
mono\-mor\-phic if and only if the base sequence is pro-\allowbreak mono\-mor\-phic.

\item The induced semidirect product sequence is semi\-stable if and only if
the base sequence is semi\-stable.

\item The induced semidirect product sequence is stable if and only if the
base sequence is stable.
\end{enumerate}

\begin{proof}
By repeated application of Lemma \ref{Lem: HNN induced triangle} we may obtain
a ladder diagram of the following type:%
\[
\begin{diagram}
G_{0}\rtimes_{\varphi_{0}}\left\langle t_{0}\right\rangle  & & \lTo^{\bar{\lambda}_{1}} & & G_{1}\rtimes_{\varphi_{1}}\left\langle
t_{1}\right\rangle  & & \lTo^{\bar{\lambda}_{2}} & &
G_{2}\rtimes_{\varphi_{2}}\left\langle t_{2}\right\rangle  & & \cdots\\
& \luTo & & \ldTo & & \luTo & & \ldTo   & \\
& & \operatorname{Im}\left(  \lambda_{1}\right)  \rtimes_{\psi_{0}}\left\langle
t_{0}\right\rangle  & & \lTo & & \operatorname{Im}\left(  \lambda
_{2}\right)  \rtimes_{\psi_{1}}\left\langle t_{1}\right\rangle  & &\lTo
& & \cdots &
\end{diagram}
\]
To prove (1), first assume that $\left\{  G_{i},\lambda_{i}\right\}  $ is
pro-\allowbreak mono\-mor\-phic. Then, after passing to a subsequence and
relabeling, we may assume that the bonding maps in the corresponding sequence
$\operatorname{Im}\left(  \lambda_{1}\right)  \leftarrow\operatorname{Im}%
\left(  \lambda_{2}\right)  \leftarrow\operatorname{Im}\left(  \lambda
_{3}\right)  \leftarrow\cdots$ are monomorphisms. Elementary properties of
semidirect products then ensure that the maps in the bottom row of the ladder
diagram are also monomorphisms. It follows that the induced semidirect product
sequence is pro-\allowbreak mono\-mor\-phic.

For the converse part of (1), notice that (even before passing to
subsequences), $\operatorname{Im}\left(  \bar{\lambda}_{i}\right)
=\operatorname{Im}\left(  \lambda_{i}\right)  \rtimes_{\psi_{0}}\left\langle
t_{i-1}\right\rangle \leq G_{i-1}\rtimes_{\varphi_{1}}\left\langle
t_{i-1}\right\rangle $, for all $i$. So by hypothesis, after passing to
subsequences and relabeling, we may assume the existence of a sequence of the
form%
\[
\operatorname{Im}\left(  \lambda_{1}\right)  \rtimes_{\psi_{0}}\left\langle
t_{0}\right\rangle \leftarrowtail\operatorname{Im}\left(  \lambda_{2}\right)
\rtimes_{\psi_{1}}\left\langle t_{1}\right\rangle \leftarrowtail
\operatorname{Im}\left(  \lambda_{3}\right)  \rtimes_{\psi_{2}}\left\langle
t_{2}\right\rangle \leftarrowtail\cdots.
\]
Then each of the above monomorphisms restricts to a monomorphism of
$\operatorname{Im}\left(  \lambda_{i+1}\right)  $ into $\operatorname{Im}%
\left(  \lambda_{i}\right)  $. It follows that $\left\{  G_{i},\lambda
_{i}\right\}  $ is pro-\allowbreak mono\-mor\-phic.

A similar strategy can be used to obtain (2), with a key ingredient being that
the homomorphisms $\operatorname{Im}\left(  \lambda_{i}\right)  \rtimes
_{\psi_{i}}\left\langle t_{i}\right\rangle \rightarrow\operatorname{Im}\left(
\lambda_{i-1}\right)  \rtimes_{\psi_{i-1}}\left\langle t_{i-1}\right\rangle $
are surjective if and only if $\operatorname{Im}\left(  \lambda_{i}\right)
\rightarrow\operatorname{Im}\left(  \lambda_{i-1}\right)  $ is surjective.

Lastly, one may obtain (3) as a consequence of (1) and (2).
\end{proof}
\end{proposition}

As noted earlier, the situation for general HNN\ extensions is more
complicated. Our best analog of Proposition
\ref{Prop: semidirect vs base sequences} is the following:

\begin{proposition}
\label{Prop: HNN vs base sequences}Assume we are given Conditions (a)-(c). Then

\begin{enumerate}
\item If the induced HNN sequence is pro-\allowbreak mono\-mor\-phic, then the
base sequence is pro-\allowbreak mono\-mor\-phic.

\item If the base sequence is semi\-stable, then the induced HNN sequence is semi\-stable.
\end{enumerate}

\begin{proof}
First we prove (2). By hypothesis, after passing to a subsequence and
relabeling, we may assume an inverse sequence of surjections:%
\[
\operatorname{Im}\left(  \lambda_{1}\right)  \twoheadleftarrow
\operatorname{Im}\left(  \lambda_{2}\right)  \twoheadleftarrow
\operatorname{Im}\left(  \lambda_{3}\right)  \twoheadleftarrow\cdots.
\]

By repeated application of Lemma \ref{Lem: HNN induced triangle}, obtain a
corresponding ladder diagram of the form:%
\[
\begin{diagram}
G_{0}\ast_{\varphi_{0}} & & \lTo^{\bar{\lambda}_{1}} & &
G_{1}\ast_{\varphi_{1}} & & \lTo^{\bar{\lambda}_{2}} & &
G_{2}\ast_{\varphi_{2}} & & \lTo^{\bar{\lambda}_{3}} & &
\cdots\\
& \luTo & & \ldTo & & \luTo & & \ldTo  & & \luTo \\
& & \operatorname{Im}\left(  \lambda_{1}\right)  \ast_{\psi_{0}} & &
\lTo & & \operatorname{Im}\left(  \lambda_{2}\right)  \ast_{\psi_{1}}
& & \lTo & & \operatorname{Im}\left(  \lambda_{3}\right)  \ast_{\psi
_{2}} & & \cdots
\end{diagram}
\]
Since each $\operatorname{Im}\left(  \lambda_{i+1}\right)  \ast_{\psi_{i}}$ is
generated by the elements of $\operatorname{Im}\left(  \lambda_{i+1}\right)  $
together with the stable letter, the homomorphisms in the bottom row are
surjective. Thus, $\left\{  G_{i}\ast_{\varphi_{i}},\bar{\lambda}_{i}\right\}
$ is semi\-stable.

In proving (1), a key difference between this and Proposition
\ref{Prop: semidirect vs base sequences} becomes important. The
\textquotedblleft up maps\textquotedblright\ in the current ladder diagram
need not be injective, so the groups in the bottom row are not just the images
of the maps in the top row (as was the case in Proposition
\ref{Prop: semidirect vs base sequences}). Instead of proceeding in the manner
of the earlier proof, we give an argument based on first principles.

Suppose the base sequence $\left\{  G_{i},\lambda_{i}\right\}  $ is not
pro-\allowbreak mono\-mor\-phic. Then for every $j_{1}$ there exists
$j_{2}>j_{1}$ such that for arbitrarily large $j_{3}>j_{2}$ there exists a
$g\in G_{j_{3}}$ such that $g\in\ker\lambda_{j_{3}j_{1}}$ but $g\notin
\ker\lambda_{j_{3}j_{2}}$. Since each $\bar{\lambda}_{i}$ acts as $\lambda
_{i}$ on $G_{i}<G_{i}\ast_{\varphi_{i}}$, one sees the same
non-pro-\allowbreak mono\-mor\-phic behavior in $\left\{  G_{i}\ast
_{\varphi_{i}},\bar{\lambda}_{i}\right\}  $.
\end{proof}
\end{proposition}

Proposition \ref{Prop: HNN vs base sequences} may be most interesting for what
is \emph{not} included. The absence of a converse for assertion (2) has
significant implications for this paper, as does the corresponding existence
of a converse in Proposition \ref{Prop: semidirect vs base sequences}. For
completeness, we provide examples covering the missing implications\ of
Proposition \ref{Prop: HNN vs base sequences}

\begin{example}
Here we describe a situation where the base sequence is stable, but the
induced HNN sequence is not pro-\allowbreak mono\-mor\-phic. Begin with the
sequence
\[
\mathbb{Z}\overset{\operatorname{id}}{\longleftarrow}\mathbb{Z}\overset
{\operatorname{id}}{\longleftarrow}\mathbb{Z}\overset{\operatorname{id}%
}{\longleftarrow}\cdots\text{.}%
\]
For the $i^{\text{th}}$ copy of $\mathbb{Z}$, let $H_{i}=H_{i}^{\prime}%
=2^{i}\mathbb{Z}$ and let $\varphi_{i}$ be the identity homomorphism. Then the
resulting HNN sequence%
\[
\mathbb{Z\ast}_{\varphi_{1}}\overset{\overline{\operatorname{id}}%
}{\longleftarrow}\mathbb{Z\ast}_{\varphi_{2}}\overset{\overline
{\operatorname{id}}}{\longleftarrow}\mathbb{Z\ast}_{\varphi_{3}}%
\overset{\overline{\operatorname{id}}}{\longleftarrow}\cdots\text{.}%
\]
consists of non-monomorphic surjections---a situation never pro-isomorphic to
a sequence of monomorphisms.
\end{example}

\begin{example}
\label{Ex: Baumslag-Solitar example}Here is a situation where a
non-semi\-stable base sequence gives rise to an induced HNN sequence that is
stable. Let $G_{0}=\left\langle a\right\rangle $ be an infinite cyclic group,
$G_{1}=\left\langle a^{2}\right\rangle $, and $\lambda$ the inclusion map. Let
$H_{0}=\left\langle a^{2}\right\rangle $, $H_{0}^{\prime}=\left\langle
a\right\rangle $, and $\varphi_{0}:\left\langle a^{2}\right\rangle
\rightarrow\left\langle a\right\rangle $ the isomorphism taking $a^{2}$ to
$a$. Similarly, let $H_{1}=\left\langle a^{4}\right\rangle $, $H_{1}^{\prime
}=\left\langle a^{2}\right\rangle $, and $\varphi_{1}:\left\langle
a^{4}\right\rangle \rightarrow\left\langle a^{2}\right\rangle $ the
isomorphism taking $a^{4}$ to $a^{2}$. Then we have presentations%
\begin{align*}
G_{0}\ast_{\varphi_{0}}  &  =\left\langle a,t\mid t^{-1}a^{2}t=a\right\rangle
\text{, and }\\
G_{1}\ast_{\varphi_{1}}  &  =\left\langle a^{2},t\mid t^{-1}a^{4}%
t=a^{2}\right\rangle \text{.}%
\end{align*}
Each is the well-known Baumslag-Solitar group $BS\left(  2,1\right)  $.
Clearly, the induced homomorphism $\overline{\lambda}:G\ast_{\varphi_{1}%
}\rightarrow G\ast_{\varphi_{0}}$ takes $t$ to $t$ and $a^{2}$ to $a^{2}$.
Applying Tietze transformations to the presentation for $G_{1}\ast
_{\varphi_{1}}$ yields%
\begin{align*}
\left\langle a^{2},t\mid t^{-1}a^{4}t=a^{2}\right\rangle  &  \rightsquigarrow
\left\langle a^{2},t,b\mid t^{-1}a^{4}t=a^{2},b=t^{-1}a^{2}t\right\rangle \\
&  \rightsquigarrow\left\langle a^{2},t,b\mid b^{2}=tbt^{-1},tbt^{-1}%
=a^{2}\right\rangle \\
&  \rightsquigarrow\left\langle t,b\mid b^{2}=tbt^{-1}\right\rangle
\end{align*}
With respect to the final presentation, $\overline{\lambda}$ takes $t$ to $t$
and $b$ to $t^{-1}a^{2}t=a$. Thus, $\overline{\lambda}$ is an isomorphism.

By applying the above observation inductively, we can begin with the
pro-\allowbreak mono\-mor\-phic but non-semi\-stable sequence%
\[
\left\langle a\right\rangle \hookleftarrow\left\langle a^{2}\right\rangle
\hookleftarrow\left\langle a^{4}\right\rangle \hookleftarrow\left\langle
a^{8}\right\rangle \hookleftarrow\cdots
\]
and end up with a corresponding HNN sequence%
\[
\left\langle a,t\mid t^{-1}a^{2}t=a\right\rangle \leftarrow\left\langle
a^{2},t\mid t^{-1}a^{4}t=a^{2}\right\rangle \leftarrow\left\langle a^{4},t\mid
t^{-1}a^{8}t=a^{4}\right\rangle \leftarrow\cdots
\]
of isomorphisms between copies of $BS\left(  2,1\right)  $.
\end{example}

We conclude this section with an observation covering a very special case of
an induced HNN sequence. It will be used in \S \ref{Sec: Z+Z}.

\begin{lemma}
\label{Lemma: induced HNN sequence-stable case}Suppose that, in addition to
conditions (a)-(c), the inverse sequences: $G_{0}\overset{\lambda_{1}%
}{\longleftarrow}G_{1}\overset{\lambda_{2}}{\longleftarrow}G_{2}%
\overset{\lambda_{3}}{\longleftarrow}\cdots$ and $H_{0}\overset{\mu_{1}%
}{\longleftarrow}H_{1}\overset{\mu_{2}}{\longleftarrow}H_{2}\overset{\mu_{3}%
}{\longleftarrow}\cdots$ (and hence $H_{0}^{\prime}\overset{\mu_{1}^{\prime}%
}{\longleftarrow}H_{1}^{\prime}\overset{\mu_{2}^{\prime}}{\longleftarrow}%
H_{2}^{\prime}\overset{\mu_{3}^{\prime}}{\longleftarrow}\cdots$) are stable,
where $\mu_{i}$ and $\mu_{i}^{\prime}$ denote appropriate restrictions of
$\lambda_{i}$ \ Then $\underline{\varphi}=\left(  \varphi_{i}\right)
_{i\geq0}$ defines an isomorphism between subgroups $\underline{H}%
=\underleftarrow{\lim}\left\{  H_{i},\mu_{i}\right\}  $ and $\underline
{H}^{\prime}=\underleftarrow{\lim}\left\{  H_{i}^{\prime},\mu_{i}^{\prime
}\right\}  $ of $\underline{G}=\underleftarrow{\lim}\left\{  G_{i},\lambda
_{i}\right\}  $; moreover, the induced HNN sequence $G_{0}\ast_{\varphi_{0}%
}\overset{\bar{\lambda}_{1}}{\longleftarrow}G_{1}\ast_{\varphi_{1}}%
\overset{\bar{\lambda}_{2}}{\longleftarrow}G_{2}\ast_{\varphi_{2}}%
\overset{\bar{\lambda}_{3}}{\longleftarrow}\cdots$ is stable and
pro-isomorphic to the constant sequence $\left\{  \underline{G}\ast
_{\underline{\varphi}}\right\}  $.

\begin{proof}
By the stability if $\left\{  G_{i},\lambda_{i}\right\}  $ we may, after
passing to a subsequence and relabeling, assume that
\[%
\begin{array}
[c]{ccccccc}%
\operatorname{Im}(\lambda_{1}) & \overset{\cong}{\longleftarrow} &
\operatorname{Im}(\lambda_{2}) & \overset{\cong}{\longleftarrow} &
\operatorname{Im}(\lambda_{3}) & \overset{\cong}{\longleftarrow} & \cdots
\end{array}
\text{.}%
\]
By passing to a further subsequence and relabeling again, we may assume that,
within each $\operatorname{Im}(\lambda_{i})$ lies the subgroup $\lambda
_{i+1}\left(  H_{i+1}\right)  $, and the restriction of $\lambda_{i+1}$ takes
$\lambda_{i+2}\left(  H_{i+2}\right)  $ isomorphically onto $\lambda
_{i+1}\left(  H_{i+1}\right)  $ for each $i$. The analogous conditions for the
primed subgroups follow automatically. Passing to the corresponding HNN
sequence for this system yields an inverse sequence of canonical isomorphisms
between groups of the form $\operatorname{Im}(\lambda_{i+1})\ast_{\left.
\psi_{i}\right\vert }$ where $\left.  \psi_{i}\right\vert $ takes
$\lambda_{i+1}\left(  H_{i+1}\right)  $ isomorphically onto $\lambda
_{i+1}\left(  H_{i+1}^{\prime}\right)  $. The projection maps are then
isomorphisms from $\underline{G}\ast_{\underline{\varphi}}$ to
$\operatorname{Im}(\lambda_{i+1})\ast_{\left.  \psi_{i}\right\vert }$.
\end{proof}
\end{lemma}

\begin{remark}
\emph{Example \ref{Ex: Baumslag-Solitar example} shows that, in general, we
cannot expect }$\underleftarrow{\lim}\left\{  G_{i}\ast_{\varphi_{i}%
},\overline{\lambda}_{i}\right\}  $\emph{ to be isomorphic to }%
$\underleftarrow{\lim}\left\{  G_{i},\lambda_{i}\right\}  \ast_{\underline
{\varphi}}$\emph{. In that particular case, }$\underleftarrow{\lim}\left\{
G_{i}\ast_{\varphi_{i}},\overline{\lambda}_{i}\right\}  $\emph{ is isomorphic
to }$BS\left(  2,1\right)  $\emph{ while }$\underleftarrow{\lim}\left\{
G_{i},\lambda_{i}\right\}  \ast_{\underline{\varphi}}\cong\mathbb{Z}$\emph{,
since }$\underleftarrow{\lim}\left\{  G_{i},\lambda_{i}\right\}  $\emph{ is
the trivial group.}
\end{remark}

\subsection{Topology at the end of a space\label{topology}}

In this paper, all spaces are assumed to be connected, locally compact,
separable, and metrizable. A space $X$ is an \emph{ANR} (\emph{absolute
neighborhood retract}) if, whenever it is embedded as a closed subset of a
metric space $Z$, some neighborhood $U$ of $X$ retracts onto $X$. It is
well-known that manifolds and locally finite CW complexes are ANRs. While many
readers will want to focus their attention on manifolds and CW complexes,
results presented here are valid for all ANRs satisfying the stated
topological assumptions. Portions of this paper, such as Proposition
\ref{Prop: fundamental proposition}, are valid for even more general spaces

A subset $N$ of a space $X$ is a \emph{neighborhood of infinity} if
$\overline{X-N}$ is compact. Standard arguments show that, when $X$ satisfies
the conditions in the previous paragraph, and $C$ is a compact subset of $X$,
then $X-C$ contains at most finitely many unbounded components, i.e., finitely
many components with noncompact closures. If $X-C$ has both bounded and
unbounded components, the situation can be simplified by letting $C^{\prime}$
consist of $C$ together with all bounded components. Then $C^{\prime}$ is
compact, and $X-C^{\prime}$ consists entirely of unbounded components.

We say that $X$ \emph{has }$k$ \emph{ends }if there exists a compactum
$C\subseteq X$ such that, for every compactum $D$ with $C\subset D$, $X-D$ has
exactly $k$ unbounded components. When $k$ exists, it is uniquely determined;
if $k$ does not exist, we say $X$ has \emph{infinitely many ends}. Thus, a
space is $0$-ended if and only if $X$ is compact, and $1$-ended if and only if
it contains arbitrarily small connected neighborhoods of infinity.

A nested sequence $N_{0}\supseteq N_{1}\supseteq N_{2}\supseteq\cdots$ of
neighborhoods of infinity, with each $N_{i}\subseteq\operatorname*{int}%
N_{i-1}$, is \emph{cofinal }if $\bigcap_{i=0}^{\infty}N_{i}=\varnothing$. Such
a sequence is easily obtained: choose an exhaustion of $X$ by compacta
$C_{0}\subseteq C_{1}\subseteq C_{2}\subseteq\cdots$, with each $C_{i-1}%
\subseteq\operatorname*{int}C_{i}$; then let $N_{i}=X-C_{i}$.

Given a nested cofinal sequence $\left\{  N_{i}\right\}  _{i=0}^{\infty}$ of
connected neighborhoods of infinity, base points $p_{i}\in N_{i}$, and paths
$r_{i}\subset N_{i}$ connecting $p_{i}$ to $p_{i+1}$, we obtain an inverse
sequence:
\begin{equation}
\pi_{1}\left(  N_{0},p_{0}\right)  \overset{\lambda_{1}}{\longleftarrow}%
\pi_{1}\left(  N_{1},p_{1}\right)  \overset{\lambda_{2}}{\longleftarrow}%
\pi_{1}\left(  N_{2},p_{2}\right)  \overset{\lambda_{3}}{\longleftarrow}%
\cdots.\medskip\label{Defn: pro-pi1}%
\end{equation}
Here, each $\lambda_{i+1}:\pi_{1}\left(  N_{i+1},p_{i+1}\right)
\rightarrow\pi_{1}\left(  N_{i},p_{i}\right)  $ is the homomorphism induced by
inclusion followed by the change of base point isomorphism determined by
$r_{i}$. The proper ray $r:[0,\infty)\rightarrow X$ obtained by piecing
together the $r_{i}$'s in the obvious manner is referred to as the \emph{base
ray }for the inverse sequence, and the pro-isomorphism class of the inverse
sequence is called the \emph{fundamental pro-group }of\emph{ }$X$ based at $r$
and is denoted $\operatorname{pro}$-$\pi_{1}\left(  X,r\right)  $. It is a
standard fact that $\operatorname{pro}$-$\pi_{1}\left(  X,r\right)  $ is
independent of the cofinal sequence of neighborhoods $\left\{  N_{i}\right\}
$ or the base points---provided those base points tend to infinity along the
ray $r$. The pro-isomorphism class is also independent of the parameterization
of $r$ and, more generally, is independent of the proper homotopy class within
$X$ of $r$. At times we take the inverse limit of $\operatorname{pro}$%
-$\pi_{1}\left(  X,r\right)  $; the result is called the \emph{\v{C}ech
fundamental group} of $X$ based at $r$ and is denoted by $\check{\pi}%
_{1}\left(  X,r\right)  $.

Clearly, if $X$ has more than one end, the \textquotedblleft choice of ends to
which $r$ points\textquotedblright\ will affect $\operatorname{pro}$-$\pi
_{1}\left(  X,r\right)  $ (and thus, $\check{\pi}_{1}\left(  X,r\right)  $).
On a more subtle note, even if $X$ has a single end, $\operatorname{pro}$%
-$\pi_{1}$ may not be independent of base ray. For the purposes of this paper,
this issue causes no problems---all concerns can be addressed by a pair of
standard propositions given below. Both can be found \ in \cite{Ge} along with
a more thorough discussion of the the non-uniqueness issue for
$\operatorname{pro}$-$\pi_{1}$.

For simplicity, the following are formulated only for the case where $X$ is one-ended.

\begin{proposition}
\label{Prop: pro-pi semistable}Let $X$ be a one-ended space and $r:[0,\infty
)\rightarrow X$ a proper ray. Then the following are equivalent:

\begin{enumerate}
\item $\operatorname{pro}$-$\pi_{1}\left(  X,r\right)  $ is semi\-stable.

\item All proper rays in $X$ are properly homotopic.
\end{enumerate}
\end{proposition}

\begin{proposition}
\label{Prop: pro-pi pro-mono}Let $X$ be a one-ended space and $r:[0,\infty
)\rightarrow X$ a proper map. Then the following are equivalent:

\begin{enumerate}
\item $\operatorname{pro}$-$\pi_{1}\left(  X,r\right)  $ is pro-\allowbreak mono\-mor\-phic.

\item There exists a compact $C\subseteq X$ such that, for every compact set
$D$ containing $C$, there exists a compact $E$ such that every loop in $X-E$
that contracts in $X-C$ contracts in $X-D$.
\end{enumerate}
\end{proposition}

A compactum $C\subseteq X$ with the property described in Proposition
\ref{Prop: pro-pi pro-mono} is called a \emph{compact core}; thus, the
proposition may be restated to say: $\operatorname{pro}$-$\pi_{1}\left(
X,r\right)  $ is pro-\allowbreak mono\-mor\-phic if and only if $X$ contains a
compact core.

Since neither condition (2) in the above two propositions involves a base ray,
it follows that (in the one-ended case) having a pro-\allowbreak
mono\-mor\-phic or semi\-stable $\operatorname{pro}$-$\pi_{1}\left(
X,r\right)  $ for some $r$ is equivalent to having pro-\allowbreak
mono\-mor\-phic or semi\-stable $\operatorname{pro}$-$\pi_{1}$ for all base
rays. For that reason, we simply say \textquotedblleft$X$ \emph{has
pro-\allowbreak mono\-mor\-phic [resp., semi\-stable] fundamental group at
infinity}\textquotedblright.

\begin{remark}
\emph{Proposition \ref{Prop: pro-pi semistable} provides even more. In the
presence of semistability, the fundamental group at infinity of a one-ended
space is well-defined up to pro-isomorphism---this is analogous to the fact
that the fundamental group of a path-connected space is well-defined up to
isomorphism. In the pro-\allowbreak mono\-mor\-phic situation this is not
always the case.}
\end{remark}

As a final preliminary, we discuss the manner in which a ladder diagram of
groups of type (\ref{basic ladder diagram}) arises in the study of fundamental
pro-groups. In doing so, we highlight the need for passage to subsequences and
address some issues related to base points and base rays. For simplicity,
consider the case of a homeomorphism $h:P\rightarrow Q$ between noncompact
spaces. For a proper ray $r:[0,\infty)\rightarrow P$, we will exhibit the
equivalence of $\operatorname{pro}$-$\pi_{1}\left(  P,r\right)  $ and
$\operatorname{pro}$-$\pi_{1}\left(  Q,h\circ r\right)  $. Rather that opting
for the most concise treatment, we give an approach that closely resembles the
situation that arises in the proof of our main theorem.

Let $\left\{  U_{i}\right\}  $ and $\left\{  V_{i}\right\}  $ be cofinal
sequences of neighborhoods of infinity in $P$ and $Q$, respectively. By
discarding an initial segment of the base ray and reparameterizing, we may
assume that $r\left(  [i,\infty)\right)  \subseteq U_{i}$ for all integers
$i\geq0$. Next we choose \textquotedblleft interlocking\textquotedblright%
\ subsequences of $\left\{  U_{i}\right\}  $ and $\left\{  V_{i}\right\}  $.
Let $V_{j_{0}}=V_{0}$, then choose $k_{0}$ sufficiently large that $h\left(
U_{k_{0}}\right)  \subseteq V_{j_{0}}$; this is possible since homeomorphisms
are proper functions. Next choose $j_{1}>j_{0}$ large enough so that
$h^{-1}\left(  V_{j_{1}}\right)  \subseteq U_{k_{0}}$; then choose
$k_{1}>k_{0}$ so that $h\left(  U_{k_{1}}\right)  \subseteq V_{j_{1}}$. By
continuing this process inductively, one obtains a pair of subsequences and a
commutative ladder diagram%
\[
\begin{diagram}
V_{j_{0}} & & \lInto & & V_{j_{1}} & &
\lInto & & V_{j_{2}} & & \lInto & & \cdots \\
& \luTo & & \ldTo & & \luTo & & \ldTo & & \luTo \\
& & U_{k_{0}} & & \lInto & & U_{k_{1}} & &
\lInto & & U_{k_{2}} & & \cdots
\end{diagram}
\]
where each \textquotedblleft up arrow\textquotedblright\ is a restriction of
$h$ and each \textquotedblleft down arrow\textquotedblright\ is a restriction
of $h^{-1}$.

We wish to apply the fundamental group functor to the above setup to get a
diagram of the form%

\[
\begin{diagram}
\pi_{1}(V_{j_{0}},q_{j_{0}}) & & \lTo & & \pi_{1}(V_{j_{1}},q_{j_{1}}) & &
\lTo & & \pi_{1}(V_{j_{2}},q_{j_{2}}) & & \lTo & & \cdots \\
& \luTo & & \ldTo & & \luTo & & \ldTo & & \luTo \\
& & \pi_{1}(U_{k_{0}},p_{k_{0}}) & & \lTo & & \pi_{1}(U_{k_{1}},p_{k_{1}}) & &
\lTo & & \pi_{1}(U_{k_{2}},p_{k_{2}}) & & \cdots
\end{diagram}
\]
This will require choices of base points and base rays. In order to define the
up\ and down\ homomorphisms resulting in a commutative diagram, those choices
must be made in a consistent manner. The first key choice was in selecting
$h\circ r$ as the base ray for $Q$. Since $r\left(  [k_{i},\infty)\right)
\subseteq U_{k_{i}}$ for all $i$, then $h\circ r\left(  [k_{i},\infty)\right)
\subseteq V_{j_{i}}$. We will not take the trouble of reparameterizing $h\circ
r$; instead for each $i\geq0$, let $p_{k_{i}}=r\left(  k_{i}\right)  $ and
$q_{j_{i}}=h\left(  p_{k_{i}}\right)  $ be the preferred base points for
$U_{k_{i}}$ and $V_{j_{i-1}}$. The horizontal bonding homomorphisms have
already been discussed; restrictions of $h$ induce the desired up
homomorphisms\ $\pi_{1}\left(  U_{k_{i}},p_{k_{i}}\right)  \rightarrow\pi
_{1}\left(  V_{j_{i}},q_{j_{i}}\right)  $. Meanwhile, restrictions of $h^{-1}$
induce homomorphism $\pi_{1}\left(  V_{j_{i}},q_{j_{i}}\right)  \rightarrow
\pi_{1}\left(  U_{k_{i-1}},p_{k_{i}}\right)  $. To get the appropriate
\textquotedblleft down homomorphism\textquotedblright\ these are composed with
change-of-base point isomorphisms $\pi_{1}\left(  U_{k_{i-1}},p_{k_{i}%
}\right)  \rightarrow\pi_{1}\left(  U_{k_{i-1}},p_{k_{i-1}}\right)  $ induced
by the paths $\left.  r\right\vert _{\left[  k_{i-1},k_{i}\right]  }$. Since
bonding homomorphisms in the bottom row were also defined using change-of-base
point isomorphisms with respect to the same set of paths, triangles of
homomorphisms with a vertex on the top row are commutative. To see that
triangles with a vertex on the bottom row also commute, it suffices to follow
the image of a single loop---the key is that the change-of-base point path in
$P$ used to define the down homomorphism is taken by $h$ to the change-of-base
point path in $Q$ used to define the horizontal bonding homomorphism.

\section{Rectangular neighborhoods of infinity and their fundamental
groups$\label{Section: Rectangular neighborhoods}$}

In this section we focus on spaces of the form $Y\times\mathbb{R}$, where $Y$
is connected and noncompact.

Since $Y\times\mathbb{R}$ contains arbitrarily large compacta of the form
$K\times\left[  u,v\right]  $, it contains arbitrarily small neighborhoods of
infinity of the form $\overline{(Y\times\mathbb{R)-}\left(  K\times\left[
u,v\right]  \right)  }$. Such a neighborhood of infinity will be called
\emph{rectangular }and denoted $R\left(  K\times\left[  u,v\right]  \right)
$. Specifically,%
\[
R\left(  K\times\left[  u,v\right]  \right)  =\overline{(Y\times
\mathbb{R)-}\left(  K\times\left[  u,v\right]  \right)  }=(\overline
{Y-K}\times\mathbb{R)\cup}\left(  Y\times(-\infty,u]\cup\lbrack v,\infty
)\right)  \text{.}%
\]
It is easy to see that rectangular neighborhoods of infinity are connected
and, by the Generalized Seifert-VanKampen Theorem (see, for example, \cite[Th.
6.2.11]{Ge}), the fundamental group of $R\left(  K\times\left[  u,v\right]
\right)  $ is isomorphic to the fundamental group of the following graph of
groups, were the vertex groups $\pi_{1}\left(  Y\right)  $ correspond to
$\pi_{1}\left(  Y\times(-\infty,u]\right)  $ and $\pi_{1}\left(
Y\times\lbrack v,\infty)\right)  $ and the edges correspond to the components
$U_{i}$ of $Y-K$ with $\Lambda_{i}=\operatorname{Im}\left(  \pi_{1}\left(
\overline{U}_{i}\right)  \rightarrow\pi_{1}\left(  Y\right)  \right)  $; the
homomorphisms of the $\Lambda_{i}$ into the vertex groups are induced by
inclusions.%
%TCIMACRO{\FRAME{dtbpF}{2.9642in}{1.2216in}{0pt}{}{}{actions-fig2.eps}%
%width 2.9642in;  height 1.2216in;  depth 0pt;  original-width 3.9156in;
%original-height 1.4885in;  cropleft "0";  croptop "1";  cropright "1";
%cropbottom "-0.0703";  filename 'actions-fig2.eps';file-properties "XNPEU";}%
%}}%
%BeginExpansion
\begin{center}
\includegraphics[
trim=0.000000in -0.104642in 0.000000in 0.000000in,
height=1.2216in,
width=2.9642in
]%
{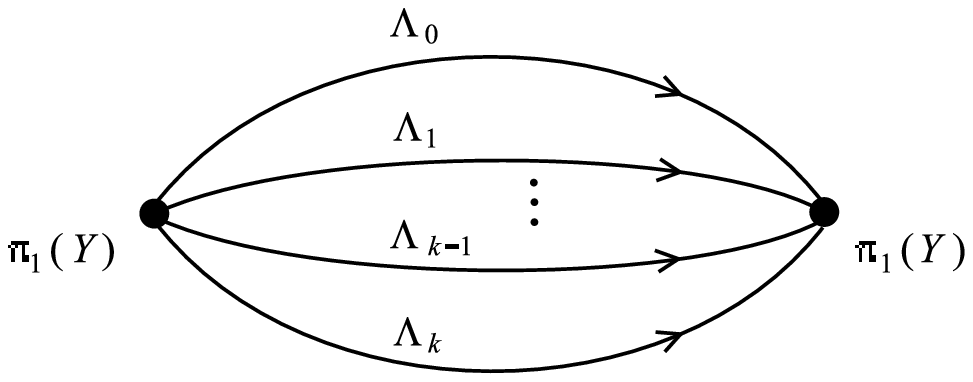}%
\end{center}
%EndExpansion

In this paper, we are particularly interested in spaces of the form $\left(
J\backslash X\right)  \times\mathbb{R}$, where $X$ is a one-ended simply
connected space and $J$ is an infinite cyclic group acting on $X$ by covering
transformations. In that case, the vertex groups will be isomorphic to
$\mathbb{Z}$ and each edge group may be viewed as $n_{i}\mathbb{Z}$ for some
integer $n_{i}$. Furthermore, by choosing $K$ so that $\left(  J\backslash
X\right)  -K$ has only unbounded components, we can ensure that each $n_{i}$
is nonzero. Indeed, if some $n_{i}$ were equal to $0$, then the corresponding
unbounded component $U_{i}$ of $\left(  J\backslash X\right)  -K$ would have
infinitely many homeomorphic preimages under the covering projection
$X\rightarrow J\backslash X$. Since $U_{i}$ has compact boundary, that would
violate the one-endedness of $X$. If $k+1$ is the number of components of
$\left(  J\backslash X\right)  -K$, the resulting graph of groups%
%TCIMACRO{\FRAME{dtbpF}{2.9337in}{1.2182in}{0pt}{}{}{actions-fig3.eps}%
%width 2.9337in;  height 1.2182in;  depth 0pt;  original-width 3.847in;
%original-height 1.4461in;  cropleft "0";  croptop "1";  cropright "1";
%cropbottom "-0.0899";  filename 'actions-fig3.eps';file-properties "XNPEU";}%
%}}%
%BeginExpansion
\begin{center}
\includegraphics[
trim=0.000000in -0.130004in 0.000000in 0.000000in,
height=1.2182in,
width=2.9337in
]%
{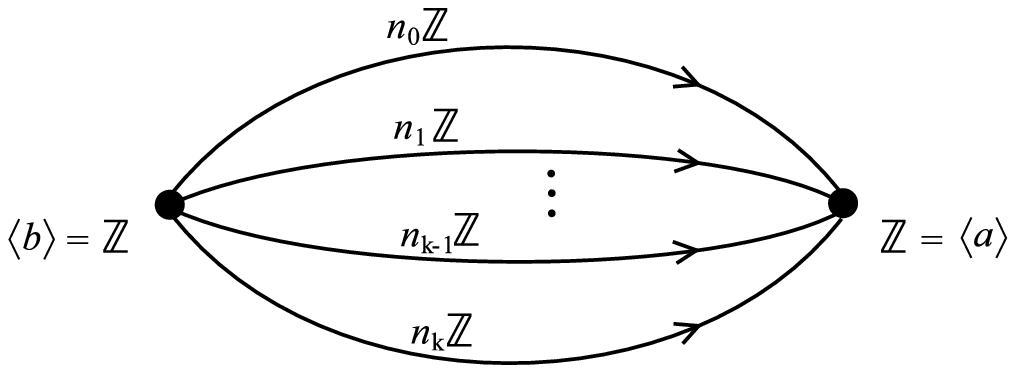}%
\end{center}
%EndExpansion
has a fundamental group $\Theta$ with presentation
\begin{equation}
\Theta=\left\langle a,b,s_{1},\cdots,s_{k}\mid a^{n_{0}}=b^{n_{0}},a^{n_{i}%
}=s_{i}b^{n_{i}}s_{i}^{-1}\text{ for }i=1,\cdots,k\right\rangle .
\label{presentation of theta}%
\end{equation}
In applying the Generalized Seifert-VanKampen Theorem to obtain this
presentation, numerous choices of base points, base paths, and a base tree
must be made. All of that is set out carefully in \cite{GM}. We provide a
brief description and refer the reader to \cite{GM} for precise details.

For a given base point $p_{0}\in J\backslash X$, the letters $a$ and $b$
represent corresponding generators of $\pi_{1}(\left(  J\backslash X\right)
\times\left\{  v+1\right\}  ,\left(  p_{0},v+1\right)  )$ and $\pi_{1}(\left(
J\backslash X\right)  \times\left\{  u-1\right\}  ,\left(  p_{0},u-1\right)
)$ (the latter connected to $\left(  p_{0},v+1\right)  $ by a base path),
respectively. If $U_{0},U_{1},\cdots U_{k}$ are the components of
$(J\backslash X)-K$, with a base point $z_{i}\in U_{i}$ for each $i$; then the
generator $s_{i}$ corresponds to a loop that passes from $(p_{0},v+1)$ to
$(z_{i},v+1)$ in $\left(  J\backslash X\right)  \times\left\{  v+1\right\}  $
then travels down the interval $z_{i}\times\left[  u-1,v+1\right]  $ to
$\left(  z_{i},u-1\right)  $ before moving within $\left(  J\backslash
X\right)  \times\left\{  u-1\right\}  $ to $\left(  z_{0},u-1\right)  $, then
up the interval $z_{0}\times\left[  u-1,v+1\right]  $ before returning to
$(p_{0},v+1)$. For consistency, all base points and all non-specified paths
are required to have first coordinate a prechosen tree $T\subseteq J\backslash
X$. For convenience, we will refer $a$ and $b$ as the \emph{primary
generators} of $\Theta$ and $s_{1},\cdots,s_{k}$ as the \emph{secondary
generators}. A key algebraic fact presented in \cite{GM} is that $\Theta$ has
a nontrivial center generated by $a^{N}$ ($=b^{N}$), where $N$ is the least
common multiple of $\left\{  n_{0},n_{1},\cdots,n_{k}\right\}  $.

Given a compactum $K^{\prime}\supseteq K$ and integers $u^{\prime
}<u<v<v^{\prime}$, we now look at the inclusion-induced homomorphism between
the fundamental groups of $R\left(  K^{\prime}\times\left[  u^{\prime
},v^{\prime}\right]  \right)  $ and $R\left(  K\times\left[  u,v\right]
\right)  $. Choose $K^{\prime}$ so that components $U_{1}^{\prime}%
,\cdots,U_{k^{\prime}}^{\prime}$ of $(J\backslash X)-K^{\prime}$ are all
unbounded, and note that whenever $U_{j}^{\prime}\subseteq U_{i}$ then $n_{i}$
divides $n_{j}^{\prime}$; By choosing base points, base paths, and base trees
that are parallel to those used for $R\left(  K\times\left[  u,v\right]
\right)  $ (in a sense made precise in \cite{GM}) we obtain a similar
presentation for the fundamental group $\Theta^{\prime}$ of $R\left(
K^{\prime}\times\left[  u^{\prime},v^{\prime}\right]  \right)  $ based at
$\left(  p_{0},v^{\prime}+1\right)  $:%
\begin{equation}
\Theta^{\prime}=\left\langle a,b,s_{1}^{\prime},\cdots,s_{k^{\prime}}^{\prime
}\mid a^{n_{0}^{\prime}}=b^{n_{0}^{\prime}},a^{n_{i}^{\prime}}=s_{i}^{\prime
}b^{n_{i}^{\prime}}s_{i}^{\prime-1}\text{ for }i=1,\cdots,k^{\prime
}\right\rangle \label{presentation of theta-prime}%
\end{equation}
Here we abuse notation slightly by again denoting the primary generators by
$a$ and $b$. This is justified by the fact that, in $R\left(  K\times\left[
u,v\right]  \right)  $ the loops representing principal generators of $\pi
_{1}(R\left(  K^{\prime}\times\left[  u^{\prime},v^{\prime}\right]  \right)
)$ are parallel copies (with respect to the product structure) of those
representing principal generators of $\pi_{1}(R\left(  K\times\left[
u,v\right]  \right)  )$. By choosing the obvious path between base points,
namely $\left\{  p_{0}\times\lbrack v+1,v^{\prime}+1]\right\}  $, the
homomorphism $\eta:\Theta^{\prime}\rightarrow\Theta$ induced by inclusion
followed by a change-of-base points isomorphism takes principal generators to
principal generators. By our notational convention, we have $\eta\left(
a\right)  =a$ and $\eta\left(  b\right)  =b$. On the secondary generators,
$\eta\left(  s_{j}^{\prime}\right)  =s_{i}$ whenever $U_{j}^{\prime}\subseteq
U_{i}$. Note that, since each $U_{i}$ contains at least one component of
$(J\backslash X)-K^{\prime}$, $\eta$ is surjective.

Now choose an exhausting sequence $K_{0}\subseteq K_{1}\subseteq
K_{2}\subseteq\cdots$ of compacta in $J\backslash X$ with each $K_{i}%
\subseteq\operatorname*{int}K_{i+1}$ and each $\left(  J\backslash X\right)
-K_{i}$ having only unbounded components. If we also choose a monotone
sequence $\{v_{i}\}$ of positive integers approaching $\infty$ and a monotone
sequence $\left\{  u_{i}\right\}  $ of negative integers approaching $-\infty$
along with base points, base paths and base trees as described above, then we
obtain an inverse sequence representing $\operatorname{pro}$-$\pi_{1}\left(
(J\backslash X)\times\mathbb{R},r\right)  $ for which all groups are of the
type described in (\ref{presentation of theta}) and
(\ref{presentation of theta-prime}) and all bonding homomorphisms are of the
type just described. In this setup, the base ray $r:[0,\infty)\rightarrow
(J\backslash X)\times\mathbb{R}$ is the one obtained by piecing together the
change of base point paths; therefore, it is a parameterization of $\left\{
p_{0}\right\}  \times\lbrack v_{0}+1,\infty)$. In a later application, it will
be convenient to use a different base ray. Like the one just described, the
new ray will pass through the base points $\left(  p_{0},v_{i}+1\right)  $ in
order and will lie entirely in $(J\backslash X)\times\lbrack v_{0},\infty)$;
but as it progresses upward with respect to the $\mathbb{R}$ component, its
$J\backslash X$ $\ $component will wind once around the loop $a$ for each unit
it moves upward. For that reason, the bonding homomorphism $\xi_{i+1}%
:\Theta_{i+1}\rightarrow\Theta_{i}$ with respect to the new ray will be the
conjugate by $a^{v_{i+1}-v_{i}}$ of the homomorphism described above. In this
scenario, the homomorphisms become less canonical (for example, except for
$a$, primary and secondary generators need not be sent to primary and
secondary generators). For similar choices of base ray $r$ (connecting the
chosen base points and lying entirely in $\left(  J\backslash X\right)
\times\lbrack v_{0},\infty)$), variations in the bonding homomorphisms
occur---in particular, conjugation by various powers of $a$. In all cases, the
properties needed later in this paper are preserved. For easy reference, we
assemble those in the following general proposition.

\begin{proposition}
\label{Prop: rectangular nbd proposition}Given $p_{0}\in J\backslash X$,
sequences $\left\{  K_{i}\right\}  $, $\left\{  u_{i}\right\}  $ and $\left\{
v_{i}\right\}  $ as described above, and base ray $r$ connecting the base
points $\left(  p_{0},v_{i}+1\right)  $ and lying entirely in $\left(
J\backslash X\right)  \times\lbrack v_{0},\infty)$, the corresponding cofinal
sequence $\left\{  R\left(  K_{i}\times\left[  u_{i},v_{i}\right]  \right)
\right\}  $ of neighborhoods of infinity in $\left(  J\backslash X\right)
\times\mathbb{R}$ gives rise to a representative of $\operatorname{pro}$%
-$\pi_{1}\left(  \left(  J\backslash X\right)  \times\mathbb{R},r\right)  $ of
the form%
\begin{equation}%
\begin{array}
[c]{ccccccc}%
\Theta_{0} & \overset{\rho_{1}}{\twoheadleftarrow} & \Theta_{1} &
\overset{\rho_{2}}{\twoheadleftarrow} & \Theta_{2} & \overset{\rho_{3}%
}{\twoheadleftarrow} & \cdots
\end{array}
\label{theta inverse sequence}%
\end{equation}
where each $\Theta_{i}$ contains a nontrivial distinguished element $a$ such that

\begin{enumerate}
\item $\rho_{i}\left(  a\right)  =a$ for all $i$, and

\item there exists a monotone sequence of positive integers $\left\{
N_{i}\right\}  $ such that, for each $i$, $a^{N_{i}}$ lies in the center of
$\Theta_{i}$.
\end{enumerate}
\end{proposition}

\begin{remark}
\emph{Note that any subsequence of \ref{theta inverse sequence} will also
satisfy the properties identified in this proposition.}
\end{remark}

\section{Mapping torus neighborhoods of infinity and their fundamental
groups\label{Section: Mapping torus neighborhoods}}

Let $X$ be a space and $j:X\rightarrow X$ a map. The mapping torus of $j$ is
the quotient space%
\begin{equation}
T_{j}\left(  X\right)  =X\times\left[  0,1\right]  /\left\{  \left(
x,1\right)  \sim\left(  j\left(  x\right)  ,0\right)  \text{ for all }x\in
X\right\}  . \label{Defn: mapping torus}%
\end{equation}
In other words, $T_{j}\left(  X\right)  $ is obtained from the mapping
cylinder of $j$ by identifying the domain and range ends of the cylinder. In
this paper, $j$ will always be a self-homeomorphism of $X$.

For all $s\in\left[  0,1\right]  $, $X\times\left\{  s\right\}
\hookrightarrow T_{j}\left(  X\right)  $ is an embedding. Amongst these
\textquotedblleft copies\textquotedblright\ of $X$ in $T_{j}\left(  X\right)
$, we give special designation to the image of $X\times\left\{  0\right\}  $.
Through this embedding, we view $X$ as a subspace of $T_{j}\left(  X\right)
$, and for $x\in X$ and $A\subseteq X$ we have $x\in T_{j}\left(  X\right)  $
and $A\subseteq T_{j}\left(  X\right)  $.

If $U$ and $V$ are subsets of $X$ with $V\subseteq U$ and $j\left(  V\right)
\subseteq U$ (thus $V\subseteq j^{-1}\left(  U\right)  $), we may define the
following subspace of $T_{j}\left(  X\right)  $:%
\[
T_{j}\left(  U,V\right)  =\left(  U\times\left\{  0\right\}  \right)
\cup\left(  V\times\lbrack0,1]\right)  /\sim
\]
where $\sim$ is the restriction of equivalence relation in
(\ref{Defn: mapping torus}). If $U$ and $V$ are neighborhoods of infinity in
$X$, then $T_{j}\left(  U,V\right)  $ is a neighborhood of infinity in
$T_{j}\left(  X\right)  $ which we refer to as a \emph{torus-like
}neighborhood of infinity. Such neighborhoods of infinity can be made
arbitrarily small by choosing $U$ to be sufficiently small in $X$. The
following proposition and its corollary give a nice description of the
fundamental group of $T_{j}\left(  U,V\right)  $ in terms of the fundamental
groups of $U$ and $V$. That description is valid only when $X$ is one-ended
with pro-\allowbreak mono\-mor\-phic fundamental group at infinity---a fact
which explains the presence of those hypotheses\ in most of our results.

\begin{proposition}
\label{key lemma}Let $X$ be connected and one-ended with pro-\allowbreak
mono\-mor\-phic fundamental group at infinity. Suppose $U$ is a connected
neighborhood of infinity such that $X-U$ contains a compact core, and
$j:X\rightarrow X$ is a homeomorphism. Then, for sufficiently small connected
neighborhoods of infinity $V\subseteq U$ and base point $p\in V:$

\begin{enumerate}
\item $j\left(  V\right)  \subseteq U$, and

\item $\ker\left(  i_{\#}\right)  =\ker\left(  \left.  j\right\vert
_{\#}\right)  $ where $i_{\#}:\pi_{1}\left(  V,p\right)  \rightarrow\pi
_{1}\left(  U,p\right)  $ and $\left.  j\right\vert _{\#}:$ $\pi_{1}\left(
V,p\right)  \rightarrow\pi_{1}\left(  U,j(p)\right)  $ are induced by
inclusion and by the restriction of $j$, respectively.
\end{enumerate}

\begin{proof}
Begin by choosing a neighborhood of infinity $U_{1}$ so small that $U_{1}$,
$j\left(  U_{1}\right)  $ and $j^{-1}\left(  U_{1}\right)  $ are all contained
in $U$. This is easy---choose $U_{1}\subseteq X-\left(  j\left(  A\right)
\cup A\cup j^{-1}\left(  A\right)  \right)  $ where $A=\overline{X-U}$.
Applying Proposition \ref{Prop: pro-pi pro-mono}, choose a smaller
neighborhood of infinity $U_{2}$ with the property that loops in $U_{2}$ which
contract in $U$ also contract in $U_{1}$. Finally, choose $V$ to be any
connected neighborhood of infinity sufficiently small that both $V$ and
$j\left(  V\right)  $ are contained in $U_{2}$.

Let $\beta$ be a loop in $V$ based at $p$. If $\beta$ contracts in $U$, then
$\beta$ also contracts in $U_{1}$; and since $j\left(  U_{1}\right)  \subseteq
U$, then $j\circ\beta$ contracts in $U$. Therefore $\ker\left(  i_{\#}\right)
\subseteq\ker(\left.  j\right\vert _{\#})$. Conversely if $j\circ\beta$
contracts in $U$ then, since $\operatorname{Im}(j\circ\beta)$ lies in $U_{2}$,
it also contracts in $U_{1}$. But $j^{-1}\left(  U_{1}\right)  \subseteq U$,
so $\beta$ also contracts in $U$. Hence, $\ker\left(  \left.  j\right\vert
_{\#}\right)  \subseteq\ker\left(  i_{\#}\right)  $.
\end{proof}
\end{proposition}

\begin{corollary}
\label{key corollary}Using the same hypotheses and notation as above, assume
$V$ is chosen to satisfy the conclusion of Proposition \ref{key lemma}. Let
$\alpha$ be a path in $U$ from $p$ to $j\left(  p\right)  $ and $\widehat
{\alpha}:\pi_{1}\left(  U,j(p)\right)  \rightarrow\pi_{1}\left(  U,p\right)  $
be the corresponding change-of-base point isomorphism. Then $\pi_{1}%
(T_{j}\left(  U,V\right)  ,p)$ is an HNN extension of $\pi_{1}(U,p)$ by $\phi$
with associated subgroups $H=i_{\#}\left(  \pi_{1}\left(  V,p\right)  \right)
\cong\pi_{1}\left(  V,p\right)  /\ker\left(  i_{\#}\right)  $ and $H^{\prime
}=\hat{\alpha}(\left.  j\right\vert _{\#}(\pi_{1}\left(  V,p\right)  )$, where
the isomorphism $\phi:H\rightarrow H^{\prime}$ is induced by $\hat{\alpha
}\circ\left.  j\right\vert _{\#}$ via an application of Proposition
\ref{key lemma}.

\begin{proof}
It is an immediate consequence of the Generalized Seifert-VanKampen Theorem
that $\pi_{1}\left(  T_{j}\left(  U,V\right)  ,p\right)  $ has the form
\[
\left\langle \pi_{1}(U,p),t\mid t^{-1}(i_{\#}\left(  g\right)  )t\left(
\hat{\alpha}(\left.  j\right\vert _{\#}(g)\right)  ^{-1}=1\text{ }%
\forall\text{ }g\in\pi_{1}\left(  V,p\right)  \right\rangle .
\]
Quotienting out by $\ker\left(  i_{\#}\right)  $ and $\ker\left(  \left.
j\right\vert _{\#}\right)  $ we obtain a presentation
\[
\left\langle \pi_{1}(U,p),t\mid t^{-1}ht\phi\left(  h\right)  ^{-1}=1\text{
}\forall\text{ }h\in H\right\rangle .
\]
Thus,
\[
\pi_{1}\left(  T_{j}\left(  U,V\right)  ,p\right)  \cong\left\langle \pi
_{1}(U,p),t\mid t^{-1}ht=\phi(h)\text{ }\forall\text{ }h\in H\right\rangle
=\pi_{1}(U,p)\ast_{\phi}\text{.}%
\]

\end{proof}
\end{corollary}

\begin{remark}
\label{Remark: HNN isomorphism}\emph{For the coming applications, it is
necessary to understand some specifics of the isomorphism }$\pi_{1}\left(
T_{j}\left(  U,V\right)  ,p\right)  \overset{\cong}{\longrightarrow}\pi
_{1}(U,p)\ast_{\phi}$\emph{ just obtained. Clearly, elements of }$\pi
_{1}(U,p)\leq\pi_{1}(U,p)\ast_{\phi}$\emph{ are paired with their natural
preimages via the inclusion }$\left(  U,p\right)  \hookrightarrow\left(
T_{j}\left(  U,V\right)  ,p\right)  $\emph{; but the remainder of the
isomorphism is not canonical---it depends on our choice of the path }$\alpha
$\emph{ from }$p$\emph{ to }$j\left(  p\right)  $\emph{. Note that a different
choice }$\beta$\emph{ may alter the subgroup }$H^{\prime}$\emph{, replacing it
with }$gH^{\prime}g^{-1}$\emph{ (and }$\phi$\emph{ by a corresponding
conjugate) where }$g\in\pi_{1}\left(  U,p\right)  $\emph{ is represented by
the loop }$\beta\ast\alpha^{-1}$\emph{. More significantly, with }$\alpha
$\emph{ as the chosen path, a loop representing the preimage of the stable
letter }$t$\emph{ is obtained by following the fiber }$\left\{  p\right\}
\times\left[  0,1\right]  $\emph{ of }$T_{j}\left(  X\right)  $\emph{ from its
initial point }$p$\emph{ to its endpoint }$j(p)$\emph{ then returning to }%
$p$\emph{ along }$\alpha^{-1}$\emph{; we denote this loop }$\tau_{p,\alpha}%
$\emph{. If }$\beta$\emph{ is used in place of }$\alpha$\emph{, the new
preimage of }$t$\emph{ is represented by }$\tau_{p,\alpha}\ast g^{-1}$\emph{.
For clarity, let }$\chi_{\alpha}$\emph{ denote the specific isomorphism}%
\[
\chi_{\alpha}:\pi_{1}\left(  T_{j}\left(  U,V\right)  ,p\right)
\overset{\cong}{\longrightarrow}\pi_{1}(U,p)\ast_{\phi}%
\]
\emph{obtained when the path }$\alpha$\emph{ is used. The key facts here are
that }$\chi_{\alpha}$\emph{ is the identity on }$\pi_{1}(U,p)$\emph{ (which,
as a result of Corollary \ref{key corollary}, may be viewed as a subgroup of
}$\pi_{1}\left(  T_{j}\left(  U,V\right)  ,p\right)  $\emph{) and takes }%
$\tau_{p,\alpha}$\emph{ to the stable letter }$t$\emph{.}
\end{remark}

\begin{proposition}
\label{Prop: Long proposition}For $X$ a connected, one-ended ANR with
pro-\allowbreak mono\-mor\-phic fundamental group at infinity, let

\begin{itemize}
\item $\left\{  U_{n}\right\}  _{n=0}^{\infty}$ be a nested cofinal sequence
of connected neighborhoods of infinity such that each pair $\left(
U_{n},U_{n+1}\right)  $ satisfies the conditions on $\left(  U,V\right)  $ in
Proposition \ref{key lemma},

\item $r:[0,\infty)\rightarrow X$ be a proper ray such that $r\left(
[n,\infty)\right)  \subseteq U_{n+1}$ for all $n$,

\item for each $n$, let $p_{n}=r\left(  n\right)  $ and $r_{n}$ be the path
from $p_{n}$ to $p_{n+1}$ obtained by restricting $r$ to $\left[
n,n+1\right]  $

\item for convenience, let $G_{n}$ denote $\pi_{1}\left(  U_{n},p_{n}\right)
$, and let $\lambda_{n+1}:G_{n+1}\rightarrow G_{n}$ be the homomorphism
$\widehat{r}_{n}\circ i_{n+1\#}$ where $i_{n+1}:\left(  U_{n+1},p_{n+1}%
\right)  \hookrightarrow\left(  U_{n},p_{n+1}\right)  $ and $\widehat{r}_{n}$
is the change-of-base points isomorphism corresponding to $r_{n}$.
\end{itemize}

\noindent Then $\operatorname{pro}$-$\pi_{1}\left(  X,r\right)  $ may be
represented by the inverse sequence%
\[
G_{0}\overset{\lambda_{1}}{\longleftarrow}G_{1}\overset{\lambda_{2}%
}{\longleftarrow}G_{2}\overset{\lambda_{3}}{\longleftarrow}\cdots.
\]
Given all of the above and a self-homeomorphism $j:X\rightarrow X$, consider
the sequence $\left\{  T_{j}\left(  U_{n},U_{n+1}\right)  \right\}
_{n=0}^{\infty}$ of connected neighborhoods of infinity in $T_{j}(X)$. For
each $n$,

\begin{itemize}
\item choose a path $\alpha_{n}$ in $U_{n}$ from $p_{n}$ to $j\left(
p_{n}\right)  $,

\item let $i_{n+1}^{\prime}:\left(  U_{n+1},p_{n}\right)  \hookrightarrow
\left(  U_{n},p_{n}\right)  $ be the inclusion map,

\item let $H_{n}=i_{n+1\#}^{\prime}\left(  \pi_{1}\left(  U_{n+1}%
,p_{n}\right)  \right)  $ and $H_{n}^{\prime}=\hat{\alpha}_{n}(\left.
j\right\vert _{\#}(\pi_{1}\left(  U_{n+1},p_{n}\right)  )$ be subgroups of
$G_{n}$ and $\phi_{n}:H_{n}\rightarrow H_{n}^{\prime}$ the isomorphism induced
by $\hat{\alpha}_{n}\circ\left.  j\right\vert _{\#}$ as promised by Corollary
\ref{key corollary},

\item let $\mu_{n}:\pi_{1}\left(  T_{j}\left(  U_{n},U_{n+1}\right)
,p_{n}\right)  \rightarrow\pi_{1}\left(  T_{j}\left(  U_{n-1},U_{n}\right)
,p_{n-1}\right)  $ be induced by inclusion followed by the change-of-base
point isomorphism corresponding to $r_{n-1}$,

\item let $\chi_{\alpha_{n}}:\pi_{1}\left(  T_{j}\left(  U_{n},U_{n+1}\right)
,p_{n}\right)  \rightarrow G_{n}\ast_{\phi_{n}}$ be the isomorphism described
in Remark \ref{Remark: HNN isomorphism}, and

\item let $\mu_{n}^{\prime}$ be the homomorphism induced by the diagram%
\begin{equation}
\begin{diagram} \pi_{1}\left( T_{j}\left( U_{n-1},U_{n}\right) ,p_{n-1}\right) & \lTo^{\mu_{n}}& \pi_{1}\left( T_{j}\left( U_{n},U_{n+1}\right) ,p_{n}\right) \\ \dTo^{\chi_{\alpha_{n-1}}} & & \dTo^{\chi_{\alpha_{n}}}\\ G_{n-1}\ast_{\phi_{n-1}} & \lDashto^{\mu_{n}^{\prime}} & G_{n}\ast_{\phi_{n}}\end{diagram} \label{HNN diagram}%
\end{equation}
Then $\operatorname{pro}$-$\pi_{1}\left(  T_{j}(X),r\right)  $ admits a
representative of the form
\[
G_{0}\ast_{\phi_{0}}\overset{\mu_{1}^{\prime}}{\longleftarrow}G_{1}\ast
_{\phi_{1}}\overset{\mu_{2}^{\prime}}{\longleftarrow}G_{2}\ast_{\phi_{2}%
}\overset{\mu_{3}^{\prime}}{\longleftarrow}\cdots
\]
where each $\mu_{i}^{\prime}$ is equal to $\lambda_{i}$ when restricted to
$G_{i}$.
\end{itemize}
\end{proposition}

To get the level of precision necessary for future arguments, we need an
additional refinement to the above representative of $\operatorname{pro}$%
-$\pi_{1}\left(  T_{j}(X),r\right)  $. That will require a hypothesis on the
homeomorphism $j:X\rightarrow X$; fortunately, that hypothesis is always
satisfied in the cases of interest in this paper.

\begin{proposition}
\label{Prop: HNN bonds are induced}In addition to the setup of Proposition
\ref{Prop: Long proposition}, assume the existence of a proper homotopy $H$ in
$X$ between the base ray $r$ and its image $j\circ r$. By reparameterizing and
rechoosing base points, if necessary, assume that $H\left(  [n,\infty
)\times\left[  0,1\right]  \right)  \subseteq U_{n}$ for all $n$. Assume also
that each $\alpha_{n}$ was chosen to be the track $H\left(  \left\{
p_{n}\right\}  \times\left[  0,1\right]  \right)  $ of the homotopy. Then each
bonding homomorphism $\mu_{n}^{\prime}$ in Proposition
\ref{Prop: Long proposition} takes the associated subgroups $H_{n}$ and
$H_{n}^{\prime}$ into $H_{n-1}$ and $H_{n-1}^{\prime}$, respectively, and
takes the stable letter $t_{n}$ to the stable letter $t_{n-1}$.

\begin{proof}
It is clear that $\mu_{n}^{\prime}$ takes $H_{n}$ into $H_{n-1}$. To see that
$\mu_{n}^{\prime}$ takes $H_{n}^{\prime}$ into $H_{n-1}^{\prime}$, we must
take note of the base paths used. An arbitrary element of $H_{n}^{\prime}$ may
be represented by a loop of the form $\alpha_{n}\ast(j\circ\gamma)\ast
\alpha_{n}^{-1}$, where $\gamma$ is a loop in $U_{n+1}$ based at $p_{n}$. Its
image in $\pi_{1}\left(  U_{n-1},p_{n-1}\right)  $ is represented by the loop
$r_{n-1}\ast\alpha_{n}\ast(j\circ\gamma)\ast\alpha_{n}^{-1}\ast r_{n-1}^{-1}$.
Since the restriction $\left.  H\right\vert _{[n-1,n]\times\lbrack0,1]}$
provides a homotopy (rel endpoints) in $U_{n}$ between the paths $r_{n-1}%
\ast\alpha_{n}$ and $\alpha_{n-1}\ast(j\circ r_{n-1})$, that loop is homotopic
in $U_{n-1}$ (rel $p_{n-1}$) to
\[
\alpha_{n-1}\ast(j\circ r_{n-1})\ast(j\circ\gamma)\ast(j\circ r_{n-1}%
)^{-1}\ast\alpha_{n-1}^{-1}=\alpha_{n-1}\ast j\circ(r_{n-1}\ast\gamma\ast
r_{n-1}^{-1})\ast\alpha_{n-1}^{-1}.
\]
The latter of the above is clearly an element of $H_{n-1}^{\prime}$.

To see that $\mu_{n}^{\prime}$ takes $t_{n}$ to $t_{n-1}$, recall from Remark
\ref{Remark: HNN isomorphism}, that a representative of $t_{n}$ in
$T_{j}\left(  U_{n},U_{n+1}\right)  $ is obtained by connecting the ends of
the mapping torus fiber $\left\{  p_{n}\right\}  \times\left[  0,1\right]  $
with the path $\alpha_{n}^{-1}$. A homotopy between this loop and the
corresponding representative of $t_{n-1}$ is apparent from the (mapped in)
annulus that may be assembled from the rectangles: $r_{n-1}\times\left[
0,1\right]  $ (a family of mapping torus fibers) and $\left.  H\right\vert
_{[n-1,n]\times\lbrack0,1]}$. Since this annulus lies in $T_{j}\left(
U_{n-1},U_{n}\right)  $, has boundary components corresponding to $t_{n-1}$
and $t_{n}$, and contains the path $r_{n-1}$ connecting $p_{n-1}$ to $p_{n}$,
the result follows.
\end{proof}
\end{proposition}

\begin{remark}
\label{Remark: culminating HNN sequence}

\begin{enumerate}
\item \emph{Proposition \ref{Prop: HNN bonds are induced} tells us that each
of the bonding homomorphism }$\mu_{n}^{\prime}:G_{n}\ast_{\phi_{n}}\rightarrow
G_{n-1}\ast_{\phi_{n-1}}$\emph{ is, in fact, the homomorphism }$\overline
{\lambda}_{n}$\emph{ induced by }$\lambda_{n}:G_{n}\rightarrow G_{n-1}$\emph{
as described in Lemma \ref{induced HNN homomorphism}. Hence, under the
hypotheses of Proposition \ref{Prop: Long proposition}, a representative }%
\[
G_{0}\ast_{\phi_{0}}\overset{\overline{\lambda}_{1}}{\longleftarrow}G_{1}%
\ast_{\phi_{1}}\overset{\overline{\lambda}_{2}}{\longleftarrow}G_{2}\ast
_{\phi_{2}}\overset{\overline{\lambda}_{3}}{\longleftarrow}\cdots
\]
\emph{of }$\operatorname{pro}$\emph{-}$\pi_{1}\left(  T_{j}\left(  X\right)
,r\right)  $\emph{ may be obtained which is, algebraically, an induced HNN
sequence (as described in \S \ref{Sec: background}) obtained from a base
sequence }%
\[
G_{0}\overset{\lambda_{1}}{\longleftarrow}G_{1}\overset{\lambda_{2}%
}{\longleftarrow}G_{2}\overset{\lambda_{3}}{\longleftarrow}\cdots
\]
\emph{representing }$\operatorname{pro}$\emph{-}$\pi_{1}\left(  X,r\right)
$\emph{.}

\item \emph{An additional property of this particular HNN sequence is that
each associated subgroup }$H_{i}\leq G_{i}\ast_{\phi_{i}}$\emph{ is precisely
the image under }$\overline{\lambda}_{i+1}$\emph{ (or, equivalently, }%
$\lambda_{i+1}$\emph{) of }$G_{i+1}$\emph{.}

\item \emph{When the homeomorphism }$j$\emph{ generates a }$\mathbb{Z}%
$\emph{-action as covering transformations on }$X$\emph{, the extra hypothesis
found in Proposition \ref{Prop: HNN bonds are induced} is easily seen to be
satisfied. Let }$p_{0}$\emph{ be a base point and }$r_{0}:\left[  0,1\right]
\rightarrow X$\emph{ be a path from }$p_{0}$\emph{ to }$j\left(  p_{0}\right)
$\emph{. For each integer }$n>0$\emph{, let }$p_{n}=j^{n}\left(  p_{0}\right)
$\emph{ and }$r_{n}=j^{n}\circ r_{0}$\emph{, a path from }$p_{n}$\emph{ to
}$p_{n+1}$\emph{. Let }$r:[0,\infty)\rightarrow X$\emph{ be the ray obtained
by gluing these paths together in the obvious manner. By proper discontinuity
of the action of }$\left\langle j\right\rangle $\emph{ on }$X$\emph{, }%
$r$\emph{ is a proper ray. Moreover, a proper homotopy between }$r$\emph{ an
}$j\circ r$\emph{ is obtained by sliding each point }$r\left(  t\right)
$\emph{ one unit forward along }$r$\emph{ to }$j\left(  r\left(  t\right)
\right)  $\emph{; hence, each }$\alpha_{n}$ \emph{is just} $r_{n}.$
\end{enumerate}
\end{remark}

We close this section with a simple special case that will be useful in
\S \ref{Sec: Z+Z}.

\begin{lemma}
\label{Lemma: torus-like nbds-stable case}Let $X$ be a connected and one-ended
with stable pro-$\pi_{1}$ and let $G=\check{\pi}_{1}\left(  X,r\right)  $,
where $r$ is a proper ray in $X$. Then any homeomorphism $j:X\rightarrow X$
induces an automorphism $\phi:G\rightarrow G$ and $T_{j}\left(  X\right)  $ is
a one-ended space with stable fundamental group at infinity and $\check{\pi
}_{1}\left(  T_{j}\left(  X\right)  ,r\right)  \cong G\rtimes_{\phi}%
\mathbb{Z}$.

\begin{proof}
By Proposition \ref{Prop: pro-pi semistable}, the base ray hypothesis of
Proposition \ref{Prop: HNN bonds are induced} is satisfied, so all of our
previous work applies. By stability, we are in position to apply Lemma
\ref{Lemma: induced HNN sequence-stable case}; moreover, in the case at hand,
the groups $\underline{G}$, $\underline{H}$ and $\underline{H}^{\prime}$ of
that lemma are all canonically isomorphic to $\check{\pi}_{1}(X,r)$. Thus, the
HNN\ extension promised there is a semidirect product isomorphic to
$G\rtimes_{\phi}\mathbb{Z}$.
\end{proof}
\end{lemma}

\section{An example}

The following example is of interest because it helps to justify all of the
effort that goes into proving Theorem \ref{Th: Stably free theorem}.

\begin{example}
Let $f:(\mathbb{S}^{1},\ast)\rightarrow(\mathbb{S}^{1},\ast)$ be degree $2$
map, and let $X^{\prime}$ be the \textquotedblleft bi-infinite mapping
telescope\textquotedblright\ of the system%
\[
\cdots\overset{f}{\longleftarrow}\mathbb{S}^{1}\overset{f}{\longleftarrow
}\mathbb{S}^{1}\overset{f}{\longleftarrow}\mathbb{S}^{1}\overset
{f}{\longleftarrow}\cdots\text{.}%
\]
Then let $X$ be the space obtained by adding a single point $p$ to compactify
the \textquotedblleft left end\textquotedblright\ of $X^{\prime}$. Since $p$
has arbitrarily small compact contractible neighborhoods in $X$ and since $X$
is locally a $2$-dimensional polyhedron at all other points of $X$, it is
clear that $X$ is a locally compact ANR. Assemble a proper base ray $r$ from
the mapping cylinder arcs corresponding to the base point $\ast$. It is easy
to see that $\operatorname{pro}$-$\pi_{1}\left(  X,r\right)  )$ is represented
by the pro-\allowbreak mono\-mor\-phic system%
\[
\mathbb{Z}\overset{\times2}{\longleftarrow}\mathbb{Z}\overset{\times
2}{\longleftarrow}\mathbb{Z}\overset{\times2}{\longleftarrow}\cdots\text{.}%
\]
which is not semi\-stable. Moreover, there is a homeomorphism $j:X\rightarrow
X$ that generates a semifree $\mathbb{Z}$-action on $X$. In particular, let
$j$ translate each mapping cylinder of $X^{\prime}$ to the right by one and
let $j\left(  p\right)  =p$.

In this situation, all of our work in
\S \ref{Section: Mapping torus neighborhoods} is valid. With a little effort,
one sees that this space $X$ and its mapping torus $T_{j}\left(  X\right)  $
provide a geometric realization of the unusual algebraic situation outlined in
Example \ref{Ex: Baumslag-Solitar example}. Since Proposition
\ref{Prop: mapping torus observation} does not apply to actions of this type,
there is no contradiction to our theorem.
\end{example}

\begin{remark}
\emph{This example provides a situation where the conclusion of Wright's
Theorem \ref{Theorem: Wright's Main Theorem} is satisfied, while the
conclusion of our Theorem \ref{Th: Stably free theorem} does not hold. Of
course, this can happen only because the $\mathbb{Z}$-action is not free.}
\end{remark}

\section{Manifold results\label{Sec: Manifold results}}

In this section, \emph{manifold} means\emph{ manifold with (possibly empty)
boundary}. A compact manifold without boundary is called \emph{closed} and a
noncompact manifold without boundary is called \emph{open}. For convenience,
all manifolds are assumed to be piecewise-linear. Analogous results may be
obtained for smooth or topological manifolds in the usual ways.

Some fundamental facts from manifold topology will be used in this section.
First, a Poincar\'{e} duality argument implies that a contractible open
$n$-manifold $M^{n}$ is always one-ended, provided $n\geq2$. A deeper result,
due to Stallings \cite{St} when $n\geq5$ and requiring the corresponding
Poincar\'{e} Conjectures \cite{MT}, \cite{Fr} for $n=3,4$, asserts that such a
manifold is homeomorphic to $\mathbb{R}^{n}$ if and only if it is simply
connected at infinity. By \emph{simply connected at infinity}, we mean that
$\operatorname*{pro}$-$\pi_{1}\left(  M^{n},r\right)  $ is pro-trivial.

We begin this section with a combined application of Proposition
\ref{Prop: mapping torus observation} and a technique found in \cite{Gu2}. The
result, which we find interesting in its own right, leads to quick new proofs
of two theorems mentioned in the introduction.

\begin{proposition}
\label{torus proposition}Let $J=\left\langle j\right\rangle \cong\mathbb{Z}$
act as covering transformations on a contractible open $n$-manifold $M^{n}$ in
an orientation preserving manner. Then the mapping torus $T_{j}\left(
M^{n}\right)  $ is homeomorphic to $\mathbb{R}^{n}\times\mathbb{S}^{1}$.

\begin{proof}
For $n\leq2$ the claim follows from classical results, so assume $n\geq3$. The
quotient space $J\backslash M^{n}$ is an orientable open $n$-manifold homotopy
equivalent to $\mathbb{S}^{1}$, so by a Poincar\'{e} duality argument (e.g.,
\cite[Prop. 3.1]{GT}), $J\backslash M^{n}$ is one-ended. By the techniques of
\S \ref{Section: Rectangular neighborhoods}, a rectangular neighborhood of
infinity $R\left(  K\times\left[  u,v\right]  \right)  $ in $\left(
J\backslash M^{n}\right)  \times\mathbb{R}$ will have a fundamental group of
the form
\[
\left\langle a,b\mid a^{m}=b^{m}\right\rangle
\]
where $a$ is a generator of $\pi_{1}\left(  \left(  J\backslash M^{n}\right)
\times\left\{  v+1\right\}  \right)  $, $b$ is the corresponding generator of
$\pi_{1}\left(  \left(  J\backslash M^{n}\right)  \times\left\{  u-1\right\}
\right)  $, and $m$ is the index of the image of $\pi_{1}\left(  \left(
J\backslash M^{n}\right)  -K\right)  $ in $\pi_{1}\left(  J\backslash
M^{n}\right)  $. Inspection of the proof in \cite{Gu2} reveals that, in this
particular situation, $\pi_{1}\left(  \left(  J\backslash M^{n}\right)
-K\right)  $ surjects onto $\pi_{1}\left(  J\backslash M^{n}\right)  $. To see
that, let $\Sigma^{1}\subseteq J\backslash M^{n}$ be a nicely embedded circle
such that inclusion is a homotopy equivalence, and let $N$ be a tubular
neighborhood of $\Sigma^{1}$. There is a homotopy of $J\backslash M^{n}$
pulling the entire space into the interior of $N$. Furthermore, there exists a
loop lying just outside $N$ that generates $\pi_{1}\left(  J\backslash
M^{n}\right)  $. The techniques found in the proof of \cite[Prop. 3.3]{Gu2}
show how that loop can be pushed into $\left(  J\backslash M^{n}\right)  -K$;
therefore $\pi_{1}\left(  \left(  J\backslash M^{n}\right)  -K\right)
\rightarrow\pi_{1}\left(  J\backslash M^{n}\right)  $ is surjective. It
follows that $\pi_{1}\left(  R\left(  K\times\left[  u,v\right]  \right)
\right)  =\allowbreak\left\langle a,b\mid a^{1}=b^{1}\right\rangle
=\allowbreak\left\langle a\right\rangle =\allowbreak\pi_{1}\left(  \left(
J\backslash M^{n}\right)  \times\mathbb{R}\right)  $.

By the above, $\left(  J\backslash M^{n}\right)  \times\mathbb{R}$ is
one-ended with stable infinite cyclic fundamental group at infinity, and
inclusion induces an isomorphism%
\[
\mathbb{\check{\pi}}_{1}\left(  \left(  J\backslash M^{n}\right)
\times\mathbb{R},r\right)  \overset{\cong}{\longrightarrow}\mathbb{\pi}%
_{1}\left(  \left(  J\backslash M^{n}\right)  \times\mathbb{R}\right)  .
\]
If we let $N^{\prime}$ be a tubular neighborhood of $\Sigma^{1}\times\left\{
0\right\}  $ in $\left(  J\backslash M^{n}\right)  \times\mathbb{R}$ and
$W^{n+1}=\left(  \left(  J\backslash M^{n}\right)  \times\mathbb{R}\right)
-int\left(  N^{\prime}\right)  $, then a standard algebraic topology argument
(use excision, the Hurewicz Theorem, and the Whitehead Theorem in the
universal cover) reveals that $\partial W^{n+1}\hookrightarrow W^{n+1}$ is a
homotopy equivalence. By Siebenmann's Open Collar Theorem, when $n+1\geq5$, or
\cite{Gu1}, when $n+1=4$, it follows that $W^{n+1}\approx\partial
W^{n+1}\times\lbrack0,\infty)$. Therefore, $\left(  J\backslash M^{n}\right)
\times\mathbb{R}$ is homeomorphic to $N^{\prime}$ with an open collar attached
to its boundary. Since $\left(  J\backslash M^{n}\right)  \times\mathbb{R}$ is
orientable, that space is homeomorphic to $\mathbb{R}^{n}\times\mathbb{S}^{1}%
$. The result now follows from Proposition
\ref{Prop: mapping torus observation}.
\end{proof}
\end{proposition}

The above result is striking when applied, for example, to the exotic
contractible open manifolds $M^{n}$ constructed in \cite{Da}. Those manifolds
have very large fundamental group at infinity, yet they admit cocompact
actions as covering transformations by groups $G$ that are finite index
subgroups of infinite Coxeter groups. Certainly, the mapping torus of
$\operatorname{id}_{M^{n}}$ is homeomorphic to $M^{n}\times\mathbb{S}^{1}$;
and thus maintains a complicated fundamental group at infinity of the form
$\operatorname*{pro}$-$\pi_{1}\left(  M^{n},r\right)  \times\mathbb{Z}$.
However, every nontrivial $g\in G$ generates a $\mathbb{Z}$-action on $M^{n}$
so by the above Proposition, whenever $g$ is orientation-preserving,
$T_{g}\left(  M^{n}\right)  $ is topologically just $\mathbb{R}^{n}%
\times\mathbb{S}^{1}$. This observation can be used to obtain a new proof of
the following theorem from \cite{GM} that was reproved in \cite{My2}.

\begin{theorem}
Let $M^{n}$ be a contractible open n-manifold not homeomorphic to
$\mathbb{R}^{n}$ and suppose a group $G$ acts as covering transformations on
$M^{n}$. Then the homomorphism $G\rightarrow\mathcal{H}\left(  M^{n}\right)  $
is injective, where $\mathcal{H(}M^{n})$ denotes the group of isotopy classes
of self-homeomorphisms of $M^{n}$. In fact, $G\rightarrow\mathcal{WH}\left(
M^{n}\right)  $ is injective, where $\mathcal{WH}\left(  M^{n}\right)  $
denotes the group of proper homotopy classes of self-homeomorphisms of $M^{n}$.

\begin{proof}
By the hypothesis, $n\geq3$. Suppose some non-trivial $g\in G$ lies in the
kernel of one of the above homomorphisms. By an algebraic topology argument,
such a $g$ must be orientation preserving. (Alternatively, that issue can be
avoided by considering $g^{2}$.) Since $g$ is non-trivial, Proposition
\ref{torus proposition} asserts that, $T_{g}\left(  M^{n}\right)  $ is
homeomorphic to $\mathbb{R}^{n}\times\mathbb{S}^{1}$. On the other hand, since
$g$ is properly homotopic to $\operatorname{id}_{M^{n}}$, that mapping torus
is proper homotopy equivalent to $M^{n}\times\mathbb{S}^{1}$. Since the
pro-fundamental group is an invariant of proper homotopy type, this can happen
only if $M^{n}$ is simply-connected at infinity, contradicting the hypothesis
that it is not homeomorphic to $\mathbb{R}^{n}$.
\end{proof}
\end{theorem}

We now employ our methods to obtain a new proof of:

\begin{theorem}
[Wright's Main Theorem]\label{Theorem: Wright's Main Theorem}Let $M^{n}$ be a
contractible open $n$-manifold with pro-\allowbreak mono\-mor\-phic
fundamental group at infinity. If $M^{n}$ admits a nontrivial action by
covering transformations, then $M^{n}\approx\mathbb{R}^{n}$.

\begin{proof}
Again we may assume $n$ is at least $3$. Any group acting by covering
transformations on a contractible manifold is torsion free; thus, $M^{n}$
admits such an action by an infinite cyclic group $J=\left\langle
j\right\rangle $. By passing to an index two subgroup if necessary, assume
that the elements of $J$ are orientation preserving homeomorphisms.

By our work in \S \ref{Section: Mapping torus neighborhoods},
$\operatorname{pro}$-$\pi_{1}\left(  T_{j}(M^{n}\right)  ,r)$ has a
representative of the form
\[
G_{0}\ast_{\phi_{0}}\overset{\overline{\lambda}_{1}}{\longleftarrow}G_{1}%
\ast_{\phi_{1}}\overset{\overline{\lambda}_{2}}{\longleftarrow}G_{2}\ast
_{\phi_{2}}\overset{\overline{\lambda}_{3}}{\longleftarrow}\cdots
\]
where
\[
G_{0}\overset{\lambda_{1}}{\longleftarrow}G_{1}\overset{\lambda_{2}%
}{\longleftarrow}G_{2}\overset{\lambda_{3}}{\longleftarrow}\cdots
\]
is a representative of $\operatorname{pro}$-$\pi_{1}(M^{n},r)$. On the other
hand, Proposition \ref{torus proposition} tells us that $\operatorname{pro}%
$-$\pi_{1}\left(  T_{j}(M^{n}\right)  ,r)$ is pro-isomorphic to the constant
sequence $\left\{  \mathbb{Z},\operatorname{id}\right\}  $. After passing to
subsequences there exists a ladder diagram of the form%
\[
\begin{diagram}
G_{k_{0}}\ast_{\phi_{k_{0}}} & & \lTo & &
G_{k_{1}}\ast_{\phi_{k_{1}}} & & \lTo & &
G_{k_{2}}\ast_{\phi_{k_{2}}} & & \lTo & &
\cdots\\
& \luTo & & \ldTo & & \luTo & & \ldTo  & & \luTo \\
& & \mathbb{Z} & & \lTo^{\operatorname{id}} & & \mathbb{Z}
& & \lTo^{\operatorname{id}}  & & \mathbb{Z} & & \lTo^{\operatorname{id}} & \cdots
\end{diagram}
\]
From this diagram, we may conclude that each homomorphism in the top row has
an infinite cyclic image. Since these homomorphisms send stable letter to
stable letter. It follows that the image of each $G_{k_{i+1}}$ in $G_{k_{i}}$
is trivial. Thus, the representative
\[%
\begin{array}
[c]{ccccccc}%
G_{k_{0}} & \overset{}{\longleftarrow} & G_{k_{1}} & \overset{}{\longleftarrow
} & G_{k_{2}} & \overset{}{\longleftarrow} & \cdots
\end{array}
\]
of $\operatorname{pro}$-$\pi_{1}(M^{n},r)$ is pro-trivial; so $M^{n}$ is
simply connected at infinity and, thus, homeomorphic to $\mathbb{R}^{n}$.
\end{proof}
\end{theorem}

\section{Proof of Theorem \ref{Th: Stably free theorem}: the algebraic
part\label{Sec: Proof of Theorem-algebraic part}}

We now begin the proof of Theorem \ref{Th: Stably free theorem}. With Theorem
\ref{Th: Wright's pro-free theorem} as a starting point, we need only prove:

\begin{proposition}
\label{Prop: fundamental proposition}Let $X$ be a one-ended simply-connected,
locally compact and locally path connected separable metric space with
pro-\allowbreak mono\-mor\-phic fundamental group at infinity. If $\mathbb{Z}$
acts by covering transformations on $X$, then $X$ has stable fundamental group
at infinity.
\end{proposition}

The proof is more intricate, but similar in spirit to our proof of Theorem
\ref{Theorem: Wright's Main Theorem}. By work done in
\S \ref{Section: Mapping torus neighborhoods}, there exists a proper ray $r$
in $X$, a representative
\begin{equation}
G_{0}\overset{\lambda_{1}}{\longleftarrow}G_{1}\overset{\lambda_{2}%
}{\longleftarrow}G_{2}\overset{\lambda_{3}}{\longleftarrow}\cdots
\label{G_i  sequence}%
\end{equation}
of $\operatorname{pro}$-$\pi_{1}\left(  X,r\right)  $, subgroups $H_{i}%
,H_{i}^{\prime}\leq G_{i}$, and isomorphisms $\phi_{i}:H_{i}\rightarrow
H_{i}^{\prime}$ with $\phi_{i-1}\circ\lambda_{i}\left(  h\right)  =\lambda
_{i}\circ\phi_{i}\left(  h\right)  $ for all $h\in H_{i}$, so that the
corresponding induced HNN\ sequence%
\begin{equation}
G_{0}\ast_{\phi_{0}}\overset{\overline{\lambda}_{1}}{\longleftarrow}G_{1}%
\ast_{\phi_{1}}\overset{\overline{\lambda}_{2}}{\longleftarrow}G_{2}\ast
_{\phi_{2}}\overset{\overline{\lambda}_{3}}{\longleftarrow}\cdots
\label{HNN sequence}%
\end{equation}
represents $\operatorname{pro}$-$\pi_{1}\left(  T_{j}\left(  X\right)
,r\right)  $, where $j:X\rightarrow X$ is a generator of the $\mathbb{Z}%
$-action on $X$.

Our work in \S \ref{Section: Rectangular neighborhoods} provides, for an
appropriately chosen base ray $r^{\prime}$, a representative of
$\operatorname{pro}$-$\allowbreak\pi_{1}\left(  \left(  J\backslash X\right)
\times\mathbb{R},r^{\prime}\right)  $ of the form
\begin{equation}%
\begin{array}
[c]{ccccccc}%
\Theta_{0} & \overset{\rho_{1}}{\twoheadleftarrow} & \Theta_{1} &
\overset{\rho_{2}}{\twoheadleftarrow} & \Theta_{2} & \overset{\rho_{3}%
}{\twoheadleftarrow} & \cdots
\end{array}
\label{theta sequence}%
\end{equation}
where the $\Theta_{i}$ and the $\rho_{i}$ are as described in that section.

To prove Proposition \ref{Prop: fundamental proposition}, we need to show that
sequence ($\ref{G_i sequence}$) is semi\-stable. By Proposition
\ref{Prop: mapping torus observation} we know that $\left(  J\backslash
X\right)  \times\mathbb{R}\approx T_{j}\left(  X\right)  $. If a homeomorphism
$h:\left(  J\backslash X\right)  \times\mathbb{R}\rightarrow T_{j}\left(
X\right)  $ can be chosen which sends an appropriately chosen base ray
$r^{\prime}$ in $\left(  J\backslash X\right)  \times\mathbb{R}$ (as described
in Proposition \ref{Prop: rectangular nbd proposition}) to an appropriate base
ray $r$ in $T_{j}\left(  X\right)  $ (as described in Proposition
\ref{Prop: HNN bonds are induced} and the remark that follows it) then the
sequences (\ref{HNN sequence}) and (\ref{theta sequence}) are pro-isomorphic.
Assuming this for the moment, a quick conclusion to our proof might be
expected as follows: First use the (explicit) semistability of
(\ref{theta sequence}) to conclude that (\ref{HNN sequence}) is semistable,
then use the semistability of (\ref{HNN sequence}) to conclude that
($\ref{G_i sequence}$) is semistable. Unfortunately, Example
\ref{Ex: Baumslag-Solitar example} shows that the second of these implications
is not automatic. Instead, we require a more delicate argument that relies on
\textquotedblleft normal forms\textquotedblright\ for HNN extensions and some
special properties of the groups and diagrams at hand.

We save for the following section, an investigation of the natural
homeomorphism between $\left(  J\backslash X\right)  \times\mathbb{R}$ and
$T_{j}\left(  X\right)  $, where it will be observed that there is no problem
with base rays---and hence, there exist subsequences of (\ref{HNN sequence})
and (\ref{theta sequence}) and a commutative diagram of the form:%

\begin{equation}
\begin{diagram} G_{n_{0}} \ast_{\phi_{n_{0}}} & & \lTo^{\bar{\lambda}_{n_{0}+1,n_{1}}} & & G_{n_{1}}\ast_{\phi_{n_{1}}} & & \lTo^{\bar{\lambda}_{n_{1}+1,n_{2}}} & & G_{n_{2}}\ast_{\phi_{n_{2}}} & & \lTo^{\bar{\lambda}_{n_{2}+1,n_{3}}}& & G_{n_{3}}\ast_{\phi_{n_{3}}} & \cdots\\ & \luTo ^{u_{0}} & & \ldTo^{d_{1}} & & \luTo ^{u_{1}} & & \ldTo^{d_{2}} & & \luTo^{u_{2}} & & \ldTo^{d_{3}} &\\ & & \Theta_{m_{0}} & & \lTo^{\rho_{m_{0}+1,m_{1}}} & & \Theta_{m_{1}} & & \lTo^{\rho_{m_{1}+1,m_{2}}}& & \Theta_{m_{2}} & & \lTo^{\rho_{m_{2}+1,m_{3}}} & & \cdots \end{diagram} \label{ladder diagram}%
\end{equation}
\noindent That investigation will also show that each\ $u_{i}$ can be arranged
to take the preferred generator $a\in\Theta_{m_{i}}$ to the stable letter
$t_{n_{i}}\in G_{n_{i}}\ast_{\phi_{_{n_{i}}}}$. For now, we assume those
arrangements have been made and proceed with the algebraic part of the proof.

By our work in \S \ref{Section: Rectangular neighborhoods}, for every $i$,
there is an integer $N_{m_{i}}>0$ for which $a^{N_{m_{i}}}$ lies in the center
of $\Theta_{m_{i}}$; therefore, $t_{n_{i}}^{N_{m_{i}}}$ is central in
$\operatorname*{Im}u_{i}=\operatorname*{Im}\bar{\lambda}_{n_{i}+1,n_{i+1}}$.
For each $i$, let $K_{i}\leq G_{n_{i}}$ be the image of $G_{n_{i+1}}$ under
$\bar{\lambda}_{n_{i}+1,n_{i+1}}$ and note that $\operatorname*{Im}%
\bar{\lambda}_{n_{i}+1,n_{i+1}}$ is precisely the subgroup of $G_{n_{i}}%
\ast_{\phi_{_{n_{i}}}}$ generated by $K_{i}$ together with the stable letter
$t_{n_{i}}$. We indicate this by writing $\operatorname*{Im}\bar{\lambda
}_{n_{i}+1,n_{i+1}}=\left\langle K_{i},t_{n_{i}}\right\rangle $. Then, by
semistability, we may assume (after passing to a further subsequence and
relabeling) that each of the bonds in the corresponding pro-equivalent inverse
sequence is surjective.%
\begin{equation}
\left\langle K_{0},t_{n_{0}}\right\rangle \twoheadleftarrow\left\langle
K_{1},t_{n_{1}}\right\rangle \twoheadleftarrow\left\langle K_{2},t_{n_{2}%
}\right\rangle \twoheadleftarrow\left\langle K_{3},t_{n_{3}}\right\rangle
\twoheadleftarrow\cdots\label{K,t sequence}%
\end{equation}

For the moment let $i$ be fixed and, to simplify notation, let $N$ denote
$N_{m_{i}}$ and $t$ denote $t_{n_{i}}$. Then $t^{N}$ is central in
$\left\langle K_{i},t\right\rangle $, so $t^{-N}K_{i}t^{N}=K_{i}$; and since
$K_{i}\leq H_{n_{i}}$, an application of the standard relators in $G_{n_{i}%
}\ast_{\phi_{_{n_{i}}}}$shows that $t^{-(N-1)}\phi_{n_{i}}(K_{i})t^{N-1}%
=K_{i}$. If $N=1$ this means $\phi_{n_{i}}$ is an automorphism of $K_{i}$;
otherwise a normal forms argument implies that $\phi_{n_{i}}(K_{i})\leq
H_{n_{i}}$ and $t^{-(N-2)}\phi_{n_{i}}^{2}(K_{i})t^{N-2}=K_{i}$. Continuing
inductively, we deduce that $\phi_{n_{i}}^{q-1}(K_{i})\leq H_{m_{i}}$ and
$t^{-(N-q)}\phi_{n_{i}}^{q}(K_{i})t^{N-q}=K_{i}$ for all $1\leq q\leq N$, with
the final observation being that $\phi_{n_{i}}^{N}(K_{i})=K_{i}$. From this it
is easy to see that conjugating by any (positive or negative) power of $t$
simply permutes the subgroups of $H_{m_{i}}$ in the collection $\left\{
\phi_{n_{i}}^{r}(K_{i})\right\}  _{r=0}^{N-1}$. And if we let $K_{i}^{\prime
}=\left\langle K_{i},\phi_{n_{i}}(K_{i}),\phi_{n_{i}}^{2}(K_{i}),\cdots
,\phi_{n_{i}}^{N-1}(K_{i})\right\rangle $, the subgroup of $H_{n_{i}}$
generated by these groups, then $\left\langle K_{i},t\right\rangle
=K_{i}^{\prime}\cdot\left\langle t\right\rangle $. Since the latter two groups
intersect trivially and $K_{i}^{\prime}$ is normal in $\left\langle
K_{i},t\right\rangle $, it follows that $\left\langle K_{i},t\right\rangle $
is a semidirect product $K_{i}^{\prime}\rtimes\left\langle t\right\rangle $.

Rewrite (\ref{K,t sequence}) as
\[
K_{0}^{\prime}\rtimes\left\langle t_{n_{0}}\right\rangle \twoheadleftarrow
K_{1}^{\prime}\rtimes\left\langle t_{n_{1}}\right\rangle \twoheadleftarrow
K_{2}^{\prime}\rtimes\left\langle t_{n_{2}}\right\rangle \twoheadleftarrow
K_{3}^{\prime}\rtimes\left\langle t_{n_{3}}\right\rangle \twoheadleftarrow
\cdots
\]
and recall that each bonding homomorphism takes $t_{n_{i+1}}$ to $t_{n_{i}}$
and $K_{i+1}^{\prime}$ into $K_{i}^{\prime}$. Such homomorphisms can be
surjective only if the $K_{i+1}^{\prime}$ surject onto $K_{i}^{\prime}$.
Moreover, since $K_{i+1}^{\prime}\leq H_{n_{i}+1}\leq G_{n_{i+1}}$ and $K_{i}$
is the image of $G_{n_{i+1}}$ under $\bar{\lambda}_{n_{i}+1,n_{i+1}}$, it
follows that $K_{i}^{\prime}=K_{i}$. This provides an inverse sequence of
surjections%
\[
K_{0}\twoheadleftarrow K_{1}\twoheadleftarrow K_{2}\twoheadleftarrow
K_{3}\twoheadleftarrow\cdots
\]
which is pro-equivalent to (\ref{G_i sequence}), as desired.

\section{Proof of Theorem \ref{Th: Stably free theorem}: topological
details\label{Sec: Proof of Theorem-topological part}}

We now lay out the topological argument needed to complete the proof of
Proposition \ref{Prop: fundamental proposition} and, thus, Theorem
\ref{Th: Stably free theorem}. Our remaining task is to prove the existence of
a ladder diagram of the form of \ref{ladder diagram} with the additional
property that $u_{i}\left(  a\right)  =t_{n_{i}}$ for each $i$.

The process described at the end of Section \ref{Sec: background} applied to a
homeomorphism $h:\left(  J\backslash X\right)  \times\mathbb{R\rightarrow
}T_{j}\left(  X\right)  $ (as promised in Proposition
\ref{Prop: mapping torus observation}), produces a ladder diagram between a
subsequences of any given representations of $\operatorname{pro}$-$\pi
_{1}\left(  \left(  J\backslash X\right)  \times\mathbb{R},r\right)  $ and
$\operatorname{pro}$-$\pi_{1}\left(  T_{j}\left(  X\right)  ,h\circ r\right)
$. Since the algebraic proof presented above used specific properties of both
(\ref{theta sequence}) and (\ref{HNN sequence}), we require a base ray $r$ in
$\left(  J\backslash X\right)  \times\mathbb{R}$ of the type specified in
Proposition \ref{Prop: rectangular nbd proposition} whose image $h\circ r$ in
$T_{j}\left(  X\right)  $ is of the type described in Remark
\ref{Remark: culminating HNN sequence}. That will be accomplished by taking a
close look at Proposition \ref{Prop: mapping torus observation}. In doing so,
it will also become clear that homomorphisms $u_{i}$ in the resulting diagram
take primary generators $a$ to stable letters $t_{n_{i}}$.

Our goal is to \textquotedblleft see\textquotedblright\ the homeomorphism
$h:(J\backslash X)\times\mathbb{R}\rightarrow T_{j}\left(  X\right)  $
promised by \ref{Prop: fundamental proposition}. Toward that end, choose a map
$F:J\backslash X\rightarrow\mathbb{S}^{1}$ that induces an isomorphism on
fundamental groups and let $\pi:X\rightarrow J\backslash X$ be the quotient
map. Then $\pi$ is a universal covering map with deck transformations
generated by $j$. Let $p:\mathbb{R\rightarrow S}^{1}$ be the universal
covering map with covering transformations generated by unit translation, and
choose $f:X\rightarrow\mathbb{R}$ to be a lift of $F\circ\pi$. For any
$A\subseteq\mathbb{R}$, let $X_{A}=f^{-1}\left(  A\right)  $. Note that, for
any unit interval $\left[  y,y+1\right]  \subseteq\mathbb{R}$, the restriction
of $\pi$ to $X_{\left[  y,y+1\right]  }$ is a quotient map that creates a copy
of $J\backslash X$ by identifying each $x\in X_{y}$ with $j\left(  x\right)
\in X_{y+1}$.

Consider the diagram,%
\begin{equation}
\begin{diagram} & & X\times\left[ 0,1\right] & & \\ & \ldTo^{q_{1}} & & \rdTo^{q_{2}} & \\ (J\backslash X)\times\mathbb{R} & & \rDashto^{h} & & T_{j}\left( X\right) \end{diagram} \label{quotient map diagram}%
\end{equation}
where $q_{1}\left(  x,u\right)  =\left(  \pi\left(  x\right)  ,f\left(
x\right)  +u\right)  $ and $q_{2}$ is the quotient map that defines
$T_{j}\left(  X\right)  $; specifically, $\left(  x,1\right)  $ is identified
with $(j\left(  x\right)  ,0)$ for each $x\in X$. Note that $q_{1}\left(
x,u\right)  =q_{1}\left(  x^{\prime},u^{\prime}\right)  $ if and only if
$\pi\left(  x\right)  =\pi\left(  x^{\prime}\right)  $ and $f\left(  x\right)
+u=f\left(  x^{\prime}\right)  +u^{\prime}$. The first of those conditions
implies that $\left\vert f\left(  x\right)  -f\left(  x^{\prime}\right)
\right\vert $ is an integer; hence, $\left\vert u-u^{\prime}\right\vert $ is
an integer. It follows that $u=u^{\prime}$ or (without loss of generality)
$u=0$ and $u^{\prime}=1$. In the first case $f\left(  x\right)  =f\left(
x^{\prime}\right)  $, implying that $\left(  x,u\right)  =\left(  x^{\prime
},u^{\prime}\right)  $. In the latter case, $f\left(  x\right)  =f\left(
x^{\prime}\right)  +1$ and, since $\pi\left(  x\right)  =\pi\left(  x^{\prime
}\right)  $, this implies that $j\left(  x^{\prime}\right)  =x$. Thus $\left(
x,0\right)  $ may be viewed as $\left(  j\left(  x^{\prime}\right)  ,0\right)
$. We conclude that, when a pair of points in $X\times\left[  0,1\right]  $ is
identified under $q_{1}$, it also identified under $q_{2}$. A similar argument
gives the converse; therefore diagram (\ref{quotient map diagram}) induces a
homeomorphism $h:(J\backslash X)\times\mathbb{R}\rightarrow T_{j}\left(
X\right)  $.

\begin{remark}
\label{Remark: elementary proof}\emph{The above does not require simple
connectivity of }$X$\emph{; the construction can be carried out more generally
by choosing }$F:J\backslash X\rightarrow S^{1}$\emph{ to induce the
epimorphism }$\pi_{1}\left(  J\backslash X,b\right)  \twoheadrightarrow
\mathbb{Z}$\emph{ with kernel equal to }$\pi_{\#}\left(  \pi_{1}\left(
X,e\right)  \right)  $\emph{, where }$e\in\pi^{-1}\left(  b\right)  $\emph{.
This provides the elementary proof of Proposition
\ref{Prop: mapping torus observation} promised in the introduction.}
\end{remark}

\begin{figure}
[ptb]
\begin{center}
\includegraphics[
height=2.0457in,
width=5.2649in
]%
{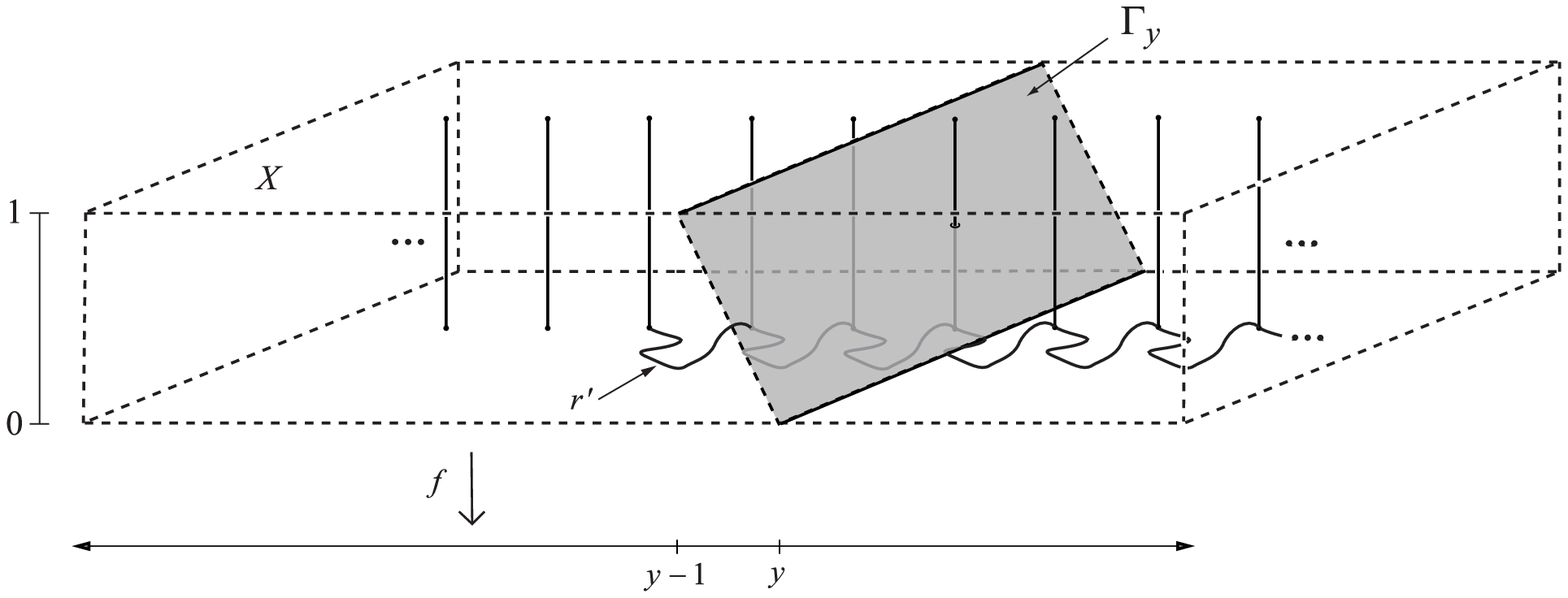}%
\caption{Preimages of $(J\backslash X)\times\left\{  r\right\}  $ and
$\left\{  b\right\}  \times\mathbb{R}$ under $q_{1}$.}%
\label{fig4}%
\end{center}
\end{figure}

Diagram (\ref{quotient map diagram}) provides a common space, $X\times\left[
0,1\right]  $, with which to compare the product structure of $(J\backslash
X)\times\mathbb{R}$ with the mapping torus structure of $T_{j}\left(
X\right)  $. The latter is easy to visualize; one simply glues the top edge
$X\times\left\{  1\right\}  $ to the bottom $X\times\left\{  0\right\}  $ via
a shift that identifies $\left(  x,1\right)  $ with $\left(  j\left(
x\right)  ,0\right)  $. To see\ the product structure of $(J\backslash
X)\times\mathbb{R}$, we look at the preimages under $q_{1}$ of factor spaces
$(J\backslash X)\times\left\{  y\right\}  $ and $\left\{  b\right\}
\times\mathbb{R}$. For fixed $y\in\mathbb{R}$,
\[
q_{1}^{-1}\left(  (J\backslash X)\times\left\{  y\right\}  \right)  =\left\{
\left(  x,u\right)  \mid u=y-f\left(  x\right)  \right\}  ,
\]
which may be viewed as the portion of the graph of the function
$(-f)+y:X\rightarrow\mathbb{R}$ lying between $u=0$ and $u=1$. Call this set
$\Gamma_{y}$ and note that it lies entirely within $X_{\left[  y-1,y\right]
}\times\lbrack0,1]$. Viewed differently, it is homeomorphic to the portion of
the graph of $f$ lying between $u=y-1$ and $u=y$ and thus is homeomorphic to
$X_{\left[  y-1,y\right]  }$. Moreover, under that homeomorphism, the
identifications made on $\Gamma_{y}$ via $q_{1}$ correspond to the
identifications made to $X_{\left[  y-1,y\right]  }$ under $\pi:X\rightarrow
J\backslash X$; both yield copies of $J\backslash X$. For fixed $b\in
J\backslash X$,%
\begin{align*}
q_{1}^{-1}\left(  (\left\{  b\right\}  \times\mathbb{R}\right)   &  =\left\{
\left(  x,u\right)  \mid\pi\left(  x\right)  =b\right\} \\
&  =\pi^{-1}\left(  b\right)  \times\left[  0,1\right] \\
&  =\left\{  j^{k}\left(  e\right)  \right\}  _{k\in\mathbb{Z}}\times\left[
0,1\right]  \text{,}%
\end{align*}
where $e\in\pi^{-1}\left(  b\right)  $. Under the quotient map $q_{1}$, the
line $\left\{  b\right\}  \times\mathbb{R}$ is assembled from $\left\{
j^{k}\left(  e\right)  \right\}  _{k\in\mathbb{Z}}\times\left[  0,1\right]  $
by identifying the top endpoint of each $j^{k}\left(  e\right)  \times\left[
0,1\right]  $ with the bottom endpoint of $j^{k+1}\left(  p\right)
\times\left[  0,1\right]  $. See Figure \ref{fig4}.%

%TCIMACRO{\FRAME{ftbpFU}{5.2649in}{2.0457in}{0pt}{\Qcb{Preimages of
%$(J\backslash X)\times\left\{  r\right\}  $ and $\left\{  b\right\}
%\times\mathbb{R}$ under $q_{1}$.}}{\Qlb{fig4}}{actions-fig4-new.eps}%
%width 5.2649in;  height 2.0457in;  depth 0pt;  original-width 7.589in;
%original-height 2.9237in;  cropleft "0";  croptop "1";  cropright "1";
%cropbottom "0";  filename 'actions-fig4-new.eps';file-properties "XNPEU";}}}%
%BeginExpansion
%EndExpansion
\bigskip

In order to obtain the desired rays $r$ and $h\circ r$ in $\left(  J\backslash
X\right)  \times\mathbb{R}$ and $T_{j}\left(  X\right)  $, we construct a
single proper ray $r^{\prime}$ in $X=X\times\left\{  0\right\}  \subseteq
X\times\left[  0,1\right]  $ and let $r=q_{1}\circ r^{\prime}$; then $h\circ
r$ is precisely $q_{2}\circ r^{\prime}$. Following the prescription found in
Remark \ref{Remark: culminating HNN sequence}, let $p_{0}\in X_{\left\{
1\right\}  }$ be a base point and $r_{0}^{\prime}:\left[  0,1\right]
\rightarrow X$ a path from $p_{0}$ to $j\left(  p_{0}\right)  \in X$; for each
$n>0$, let $p_{n}=j^{n}\left(  p_{0}\right)  $ and $r_{n}^{\prime}=j^{n}\circ
r_{0}^{\prime}$ (a path from $p_{n}$ to $p_{n+1}$). Obtain $r^{\prime
}:[0,\infty)\rightarrow X$ by gluing these paths together in the obvious
manner. See Figure \ref{fig4}. By choosing the initial path $r_{0}^{\prime}$
to lie in $X_{(0,\infty)}$, it will follow that $r^{\prime}\left(
[n,\infty)\right)  \subseteq X_{(n,\infty)}$ for all integers $n\geq0$. Since
$q_{1}\left(  X_{(n,\infty)}\times\left\{  0\right\}  \right)  \subseteq
\left(  J\backslash X\right)  \times\left(  n,\infty\right)  $, it follows
that $r\left(  [n,\infty)\right)  \subseteq\left(  J\backslash X\right)
\times\left(  n,\infty\right)  $ for all $n$; so $r$ is in accordance with
Proposition \ref{Prop: rectangular nbd proposition}. By construction,
$q_{2}\circ r^{\prime}$ fits the specifications of Remark
\ref{Remark: culminating HNN sequence}; so the base ray issue is
resolved---there is a ladder diagram of type (\ref{ladder diagram}).

Lastly we observe that each $u_{i}$ in (\ref{ladder diagram}) takes the
primary generator $a\in\Theta_{m_{i}}$ to the stable letter $t_{n_{i}}\in
G_{n_{i}}\ast_{\phi_{n_{i}}}$. Following the discussion found in Remarks
\ref{Remark: HNN isomorphism} and \ref{Remark: culminating HNN sequence},
$t_{n_{i}}$ may be represented by the loop in $T_{j}\left(  U_{n_{i}}%
,U_{n_{i}+1}\right)  $ which is a concatenation of the $q_{2}$-image of the
interval $p_{n_{i}+1}\times\left[  0,1\right]  $ with the path $\left(
\left.  r^{\prime}\right\vert _{[n_{i+1},n_{i+1}+1]}\right)  ^{-1}$ in
$X\times\left\{  0\right\}  $. The $q_{1}$-images of these paths lie in
$\left(  J\backslash X\right)  \times\lbrack n_{i},\infty)$ where the
resulting loop is easily seen to generate $\pi_{1}\left(  \left(  J\backslash
X\right)  \times\lbrack n_{i},\infty)\right)  $. This is the primary generator
$a$ of a rectangular neighborhood of infinity of the form $R\left(  K_{n_{i}%
}\times\left[  -n_{i},n_{i}\right]  \right)  $. Thus, $h$ takes this
representation of $a$ to $t_{n_{i}}$. Since, each $u_{i}$ is induced by a
restriction of the homeomorphism $h$, and since we have arranged that all
bonding homomorphisms in (\ref{theta sequence}) and (\ref{HNN sequence}) take
primary generator to primary generator and and stable letter to stable letter,
it follows that each $u_{i}$ takes $a$ to $t_{n_{i}}$.

\section{Actions by $\mathbb{Z\oplus Z}\label{Sec: Z+Z}$}

In this section we prove Theorem \ref{Th: Z+Z theorem}. Cocompact actions of
$\mathbb{Z\oplus Z}$ are well-understood, so we discuss only the non-cocompact
case. We will prove the following:

\begin{theorem}
Let $X$ be a one-ended, simply connected, locally compact ANR with
pro-\allowbreak mono\-mor\-phic fundamental group at infinity. If $X$ admits a
$\mathbb{Z\oplus Z}$ action by covering transformations that is not cocompact,
then $X$ is simply connected at infinity.

\begin{proof}
Assume $G\cong\mathbb{Z\oplus Z}$ acts by covering transformations on $X$.
Write $G=\left\langle j_{1}\right\rangle \oplus\left\langle j_{2}\right\rangle
$ where each $j_{i}$ is a self-homeomorphism of $X$. By Proposition
\ref{Prop: mapping torus observation}, $\left(  \left\langle j_{1}%
\right\rangle \backslash X\right)  \times\mathbb{R}\approx T_{j_{1}}\left(
X\right)  $; moreover, we will see that $j_{2}$ induces natural
self-homeomorphisms of each of these spaces. This is largely due to the fact
that $j_{2}j_{1}=j_{1}j_{2}$.

Let $\sim$ be the equivalence relation on $X$ induced by the action of
$\left\langle j_{1}\right\rangle $ and let $\left[  x\right]  $ denote a
corresponding equivalence class. Then $\breve{j}_{2}:\left\langle
j_{1}\right\rangle \backslash X\mathbb{\rightarrow}\left\langle j_{1}%
\right\rangle \backslash X$ defined by $\breve{j}_{2}\left(  \left[  x\right]
\right)  =[j_{2}\left(  x\right)  ]$ is a well-defined function. Indeed, if
$x\sim y$ then $y=j_{1}^{k}\left(  x\right)  $ for some integer $k$. Then
$j_{2}\left(  y\right)  =j_{2}j_{1}^{k}\left(  x\right)  =j_{1}^{k}%
j_{2}\left(  x\right)  $, and the last of these terms is equivalent to
$j_{2}\left(  x\right)  $ by definition. Continuity of $\breve{j}_{2}$ is
clear; moreover, $j_{2}^{-1}$ induces a continuous inverse for $\breve{j}_{2}$
in an analogous manner. Thus, $\breve{j}_{2}$ is a homeomorphism. Let
$\overline{j}_{2}=\left(  \breve{j}_{2},\operatorname{id}_{\mathbb{R}}\right)
:\left(  \left\langle j_{1}\right\rangle \backslash X\right)  \times
\mathbb{R\rightarrow}\left(  \left\langle j_{1}\right\rangle \backslash
X\right)  \times\mathbb{R}$.

The desired self-homeomorphism of $T_{j_{1}}\left(  X\right)  =X\times\left[
0,1\right]  /\left\{  \left(  x,1\right)  \sim\left(  j_{1}\left(  x\right)
,0\right)  \right\}  $ is obtained by letting $j_{2}$ act on each slice
$X\times\left\{  t\right\}  $. Since $j_{2}\left(  j_{1}\left(  x\right)
\right)  =j_{1}\left(  j_{2}\left(  x\right)  \right)  $, applying $j_{2}$ to
a pair of equivalent points $\left(  x,1\right)  $ and $\left(  j_{1}\left(
x\right)  ,0\right)  $ yields a pair of equivalent points $\left(
j_{2}(x),1\right)  $ and $\left(  j_{1}j_{2}\left(  x\right)  ,0\right)  $.
Let $\underline{j}_{2}:T_{j_{1}}\left(  X\right)  \rightarrow T_{j_{1}}\left(
X\right)  $ be the resulting homeomorphism.

A quick check of diagram \ref{quotient map diagram} shows that the induced
homeomorphism $h$ is equivariant with respect to the $\mathbb{Z}$-actions
induced by $\overline{j}_{2}$ and $\underline{j}_{2}$. So, by a second
application of Proposition \ref{Prop: mapping torus observation}, we get%
\begin{equation}
\left\langle \overline{j}_{2}\right\rangle \backslash\left(  (\left\langle
j_{1}\right\rangle \backslash X)\times\mathbb{R}\right)  \times\mathbb{R}%
\approx T_{\underline{j}_{2}}\left(  T_{j_{1}}\left(  X\right)  \right)  .
\label{Z+Z identity}%
\end{equation}
We are now prepared to employ our standard strategy; in particular, we will
use identity (\ref{Z+Z identity}) to obtain a pair of inverse sequences
representing pro-fundamental groups of these spaces. Comparison of those
sequences will reveal the desired conclusion.

The left-hand side of (\ref{Z+Z identity}) is easily seen to be homeomorphic
to
\[
(\left(  \left\langle j_{1}\right\rangle \oplus\left\langle j_{2}\right\rangle
\right)  \backslash X)\times\mathbb{R}^{2}.
\]
The space $\left(  \left\langle j_{1}\right\rangle \oplus\left\langle
j_{2}\right\rangle \right)  \backslash X$ has fundamental group isomorphic to
$\left\langle j_{1}\right\rangle \oplus\left\langle j_{2}\right\rangle $. It
is a standard fact that, given a path connected noncompact space $Y$, the
\textquotedblleft doubly stabilized\textquotedblright\ product $Y\times
\mathbb{R}^{2}$, has a stable fundamental group at infinity isomorphic to
$\pi_{1}\left(  Y\right)  $. To see this, first note that the noncompactness
of $Y$ implies that $Y\times\mathbb{R}$ is one-ended; moreover, the
fundamental group of each neighborhood of that end surjects onto $\pi
_{1}\left(  Y\times\mathbb{R}\right)  $. So the techniques of
\S \ref{Section: Rectangular neighborhoods}, applied to $\left(
Y\times\mathbb{R}\right)  \times\mathbb{R}$, show that this space has
arbitrarily small rectangular neighborhoods of infinity with fundamental group
described by a graph of groups with just two vertices and one edge---all
labeled by $\pi_{1}\left(  Y\right)  $.

To understand the fundamental group at infinity for $T_{\underline{j}_{2}%
}\left(  T_{j_{1}}\left(  X\right)  \right)  $, first use Theorem
\ref{Th: Stably free theorem} to deduce that $X$ has a stable finitely
generated free fundamental group at infinity. If $\check{\pi}_{1}\left(
X,r\right)  =F$, where $F$ is a finitely generated free group, then by Lemma
\ref{Lemma: torus-like nbds-stable case}, $T_{j_{1}}\left(  X\right)  $ has
stable pro-$\pi_{1}$ isomorphic to a semidirect product $F\rtimes\mathbb{Z}$.
Applying Lemma \ref{Lemma: torus-like nbds-stable case} a second time, we see
that $T_{\underline{j}_{2}}\left(  T_{j_{1}}\left(  X\right)  \right)  $ has a
stable fundamental group at infinity of the form $\left(  F\rtimes
\mathbb{Z}\right)  \rtimes\mathbb{Z}$.

Combining the above observations, we have $\left(  F\rtimes\mathbb{Z}\right)
\rtimes\mathbb{Z\cong}\left\langle j_{1}\right\rangle \oplus\left\langle
j_{2}\right\rangle \cong\mathbb{Z\oplus Z}$. When a semidirect product is
abelian, the factor groups must both be abelian and the product an ordinary
direct product. It follows that $F\rtimes\mathbb{Z=}F\times\mathbb{Z}$ and
$\left(  F\rtimes\mathbb{Z}\right)  \rtimes\mathbb{Z\allowbreak=\allowbreak
}\left(  F\times\mathbb{Z}\right)  \times\mathbb{Z}$. The latter group is
isomorphic to $\mathbb{Z\oplus Z}$ if and only if $F$ is the trivial group.
Hence, $X$ is simply connected at infinity.
\end{proof}
\end{theorem}

\end{document}